# A HIERARCHICAL PRECONDITIONER FOR WAVE PROBLEMS IN QUASILINEAR COMPLEXITY*

BORIS BONEV† AND JAN S. HESTHAVEN†

**Abstract.** The paper introduces a novel, hierarchical preconditioner based on nested dissection and hierarchical matrix compression. The preconditioner is intended for continuous and discontinuous Galerkin formulations of elliptic problems. We exploit the property that Schur complements arising in such problems can be well approximated by hierarchical matrices. An approximate factorization can be computed matrix-free and in a (quasi-)linear number of operations. The nested dissection is specifically designed to aid the factorization process using hierarchical matrices. We demonstrate the viability of the preconditioner on a range of 2D problems, including the Helmholtz equation and the elastic wave equation. Throughout all tests, including wave phenomena with high wavenumbers, the generalized minimal residual method (GMRES) with the proposed preconditioner converges in a very low number of iterations. We demonstrate that this is due to the hierarchical nature of our approach which makes the high wavenumber limit manageable.

## 1. Introduction.

Many engineering problems require the solution of increasingly large linear systems of equations

$$Ax = b, \tag{1.1}$$

where $A \in \mathbb{R}^{n\times n}$ is an invertible matrix and $x, b \in \mathbb{R}^n$. Typically, $A$ is sparse and arises from a finite element type discretizations of a partial differential equation (PDE) on complex domains involving a prohibitively large number of degrees of freedom. This makes the application of direct solvers intractable, as fill-in often becomes excessive and memory requirements cannot be met [10]. Moreover, the computational work required to solve such problems is too large and one would seek algorithms with linear-complexity. Consequently, other methods, such as multigrid methods or iterative methods [39] have been popular, as they promise to solve the problems in an amount of time that scales linearly with $n$. Krylov subspace methods seek to find the solution by projecting onto a subspace defined by powers of the original matrix. However, matrices arising from such problems often have unfavorable spectral propoerties, and may require many iterations to converge to an acceptable solution. Preconditioners are therefore often required to keep the number of iterations to a minimum, as it may otherwise become a bottleneck. This is especially true for time-dependent problems, where many solutions have to be computed. A preconditioner $P^{-1} \in \mathbb{R}^{n\times n}$ transforms the linear system of equations into the equivalent problem

$$P^{-1}Ax = P^{-1}b, \tag{1.2}$$

where $P^{-1}A$ is favorable spectral properties and therefore better suited for the solution using Krylov subspace methods, such as the generalized minimal residual method

*Submitted to the editors September 9, 2020.

**Funding:** This work was supported in part by the Swiss National Science Foundation under grant no. 513966.

†Department of Mathematics, EPFL, Lausanne, Switzerland (boris.bonev@epfl.ch, jan.hesthaven@epfl.ch).

(GMRES). However, preconditioners are often problem-specific and can become expensive to compute and/or apply.

Efficient preconditioning techniques have been difficult to develop for models of wave phenoma such as the Helmholtz problem [13, 15, 18]. These problems arise in acoustics, geophysics [30], aerospace [37] and other domains in science and engineering. Consequently, various methods have been proposed to precondition such problems. Examples include multigrid solvers based on a shifted Laplacian, analytic incomplete *LU* (AILU) and domain decomposition methods based on alternating Schwartz methods.

Shifted Laplacian methods for the Helmholtz problem [14, 18] modify the operator by adding complex entries to the diagonal of the matrix and solving the resulting problem using a multigrid approach. These preconditioners are easy to compute, however their performance typically deteriorates as the wavenumber increases as the operator tends to deviate increasingly in the large wavenumber limit. Moreover, it is not entirely clear how an appropriate shift should be chosen to guarantee convergence.

Other approaches include domain decomposition methods such as the classical Schwarz method. These methods consider subdomains and solve the original problem on these subdomains. The solutions are then connected via boundary conditions on the respective interfaces, which requires the solution of a coarse problem to coordinate the solution on the interfaces. For waveproblems, this becomes problematic as the operator can be resonant/singular on a subdomain, while being non-singular on the entirety of the domain [18]. Consequently, subdomains have to be sufficiently small, so that this can be avoided. However, the coarse problem then becomes inefficient, especially if the wavenumber is high, as too many small subdomains have to be considered.

Lately, other types of preconditioners such as the analytic incomplete *LU* factorization [16, 17] and sweeping preconditioners [11] have become more popular. The latter is based on a block factorization that eliminates degrees of freedom layer by layer, while compressing intermediate matrices in the hierarchical matrix framework. These methods have been particularly promising as their performance does not deteriorate with increasing wavenumbers. However they do not generalize to other wave phenomena and the applicability to complex domains and general grids remains to be demonstrated. Moreover, while these preconditioners can be constructed and applied in linear complexity, it is not straightforward to parallelize them.

This brings us back to direct methods, where rank-structured matrices have been increasingly popular to accelerate their construction [44, 42, 1, 20, 19, 26]. These include block low-rank, $\mathcal{H}$, $\mathcal{H}^2$, hiearchically off-diagonal low-rank (HODLR) and hierarchically semi-separable (HSS) matrix formats. Hierarchical matrix formats allow standard matrix arithmetic (addition, multiplication, inversion,...) to be executed in quasilinear time with respect to $n$ with some dependency on the off-diagonal rank $k$. On the other hand, these formats sacrifice flexibility to attain this efficiency. As a result, permutations and accessing individual entries are non-trivial. Moreover compression is costly unless these matrices arise from special problems such as the fast multipole method (FMM) or boundary element methods (BEM) [6].

In the context of direct solvers, there has been a substantial effort to apply hierarchical matrices to accelerate the factorization process. This is based on the idea that Schur complements arising in the factorization exhibit a hierarchical low-rank structure, which can be exploited by the above-mentioned formats [4, 2, 3, 12]. A core question here is how intermediate matrices arising in the factorization process can be obtained efficiently in a hierarchical format to achieve the desired speed-up.

One possibility is to perform the entire factorization process using hierarchical matrix arithmetic. This has been pursued in [44], where intermediate matrices are directly formed via the structured version of the so-called extend-add operation [10]. The main drawback of this approach is that this process is very complicated and only works for certain 2D mesh structures. Thus, simplified versions are proposed [40], which keep intermediate matrices dense at the cost of higher computational complexity and memory requirements.

The alternative is to use randomized techniques [34] and compress intermediate matrices by computing their action on a random Gaussian matrix [40, 20, 19]. In the case of [40, 20], the structured extend-add process is replaced by the 'skinny' extend-add operation which only operates on matrix-vector products required by the randomized compression. While [19] similarly avoids the complicated extend-add operation by using randomized compression, we were not able to replicate the performance of the algorithm. We believe there is an oversight concerning the necessity of accessing individual entries for the compression algorithm.

We propose a novel method based on the nested dissection of the computational domain and the compression of Schur complements using HSS matrices. This method avoids the formation and factorization of dense intermediate matrices by using a combination of HSS arithmetic and randomized sampling. As the goal is preconditioning, we can make some assumptions with respect to the nested dissection, which simplify the HSS arithmetic in the factorization process. More precisely, a nested dissection can be constructed, which completely decouples the interior of nodes from all other nodes on the same level. As a consequence, the support of Schur complements arising in one level of the elimination hierarchy becomes disjoint and we can form and factor intermediate block matrices (frontal matrices) without the need for a complex HSS extend-add operation or an extra sampling step. After factorization, Schur complements can be formed using the randomized algorithm [34], and a block structure can be chosen which allows access to submatrices needed in the factorization step of the parent node.

Consequently, we can construct the factorization in quasilinear time and apply the preconditioner in linear time. At the same time, both construction and application are matrix-free and thus memory-efficient. This makes the approach general-purpose as it is based on hierarchical block-factorization and the assumption of compressible Schur complements, and makes it applicable to general elliptic problems and wave problems. To demonstrate this, we apply our method to various problems including finite element discretizations of the Helmholtz equation and the frequency-domain elastic wave equations. These include traditionally difficult problems involving heterogeneous material parameters or non-standard geometries. Moreover, we demonstrate that the preconditioned GMRES converges in just a few iterations. This remains true even when wavenumbers become large. We demonstrate that is caused by our hierarchical approach, which allows us to adapt the hierarchy of the preconditioner to match the high wavenumbers.

The paper is organized as follows: In section 2 we discuss notation and introduce hierarchical matrices with a focus on HSS matrices. The construction of the preconditioner is discussed in section 3. Numerical results are presented in section 4 and conclusions follow in section 5.

**2. Preliminaries.** We introduce underlying concepts for hierarchical matrices. Our notation follows [22] and in particular, we often deal with submatrices of the matrix $A$, for which we will use the MATLAB style notation $A(I, J)$, where $I$ and $J$ are index sets. If we take the entire index-set, this is represented by the symbol ':', as in $A(:, J)$, where we have extracted entire columns of $A$.

Our treatment of HODLR and HSS matrices closely follows [27, 34, 21]. The reader is referred to these works for a detailed introduction to the topic. Moreover, much of the algorithms related to HSS matrices are readily available in [5].

**2.1. Hierarchical matrices.** In this work we focus on HSS (hierarchically semi-separable) matrices, which are a special case of HODLR (hierarchically off-diagonal low-rank) matrices. These matrices have a hierarchical structure, where certain blocks in this structure are low-rank. Given a block-partitioned matrix $A \in \mathbb{R}^{m\times n}$,

$$A = \begin{bmatrix} A_{11} & A_{12} \\ A_{21} & A_{22} \end{bmatrix},$$

the off-diagonal blocks $A_{12}$ and $A_{21}$ are assumed to be representable as low-rank matrices, whereas the diagonal blocks $A_{11}$ and $A_{22}$ recursively possess the same structure as $A$. This proceeds up to some leaf level, where the diagonal blocks are full rank. In the following, we will follow the treatment of [27] and formalize the concept of HSS matrices.

Let us introduce a hierarchical partitioning of the index set $I = \{1, 2, \ldots, m\}$:

DEFINITION 2.1 (Cluster tree). *Let $I = \{1, 2, \ldots, m\}$ be an index set and let $\mathcal{T}$ be a binary tree of maximal depth $L$, where all of its nodes $I_i^l \in \mathcal{T}$ are subsets of $I$. The superindex $l$ in $I_i^l$ denotes the level, starting from the top-level $l = 0$, and the subscript $i$ enumerates the node within that level. We call $\mathcal{T}$ a cluster tree if*
- *the root node $I^0$ is the entire index set $I$,*
- *nodes at each level $l$ are disjoint: $\forall i \neq j : I_i^l \cap I_j^l = \emptyset$,*
- *every non-leaf node $I_i^l$ with children $I_{c_1}^{l+1}, I_{c_2}^{l+1}$ is the union of its children $I_i^l = I_{c_1}^{l+1} \cup I_{c_2}^{l+1}$.*

Definition 2.1 implies that each level $I_1^l, I_2^l, \ldots$ of the cluster tree forms a partition of the entire index set $I$, assuming that the tree is balanced. Cluster trees are not necessarily balanced, which allows trees to be locally refined. Let us consider two cluster trees, $\mathcal{T}_r$ and $\mathcal{T}_c$, one for the row indices $I = \{1, 2, \ldots, m\}$ and one for the column indices $J = \{1, 2, \ldots, n\}$. We assume that the clusters are contigious, i.e. indices appear in incremental order when the tree is traversed in level-order. Such a clustering can be obtained by simply permuting the matrix accordingly. If both cluster trees have the same tree structure (not necessarily the same nodes), they induce a HODLR partitioning of the matrix $A$. Thus, we can take sibling nodes $i$ and $j$, which correspond to off-diagonal blocks $A(I_i^l, J_j^l)$ in $A$. As we expect these blocks to be low-rank, there is no need to partition them any further. If we choose a diagonal block with $i = j$, either both $I_i^l$ and $J_i^l$ are leaf nodes or both $I_i^l$ and $J_i^l$ are branch nodes. In the latter case, $A(I_i^l, J_i^l)$ is recursively partitioned into a two-by-two block matrix by the children of $I_i^l$ and $J_i^l$ until we reach the leaf level. Thus, we obtain the desired recursive structure of HODLR matrices:

DEFINITION 2.2 (HODLR matrix). *Consider a matrix $A \in \mathbb{R}^{m\times n}$ with row indices $I$, column indices $J$ and corresponding row- and column-cluster trees $\mathcal{T}_r$, $\mathcal{T}_c$ with identical tree structures. Given a small positive integer $k$, $A$ is said to be a HODLR matrix, with respect to the cluster trees $\mathcal{T}_r$, $\mathcal{T}_c$, if*

$$\forall I_i^l \in \mathcal{T}_r, J_j^l \in \mathcal{T}_c, \text{ such that } i \text{ and } j \text{ are siblings}: \quad \operatorname{rank} A(I_i^l, J_j^l) \leq k.$$

*The minimum $k$ for which this is true is called the HODLR rank of $A$.*

Typically, we consider matrices, for which $k$ is much smaller than the matrix dimensions $m$ and $n$, i.e. $k \ll m, n$. This leads to an improved storage complexity of $\mathcal{O}(kn \log n)$, if $m \approx n$ and if the depth of the hierarchy is chosen to be $\log n$. Other operations, such as matrix-vector multiplications are also accelerated. An overview can be found in [36]. It is worth noting, that for the same matrix $A$, different cluster trees can often result in very different HODLR ranks. Thus, we must put some thought into the generation of $\mathcal{T}_r$ and $\mathcal{T}_c$. A simple way of forming cluster trees is by recursive bisection of the initial index set $I$. This process is repeated until the size of the index set is below a certain threshold. We call

$$\beta = \max_{I_i^l \in \mathcal{T}_r \text{ or } I_i^l \in \mathcal{T}_c} |I_i^l| \tag{2.1}$$

the block-size of the partitioning. Ideally, we want to choose a block-size (2.1), that is large enough to guarantee efficient low-rank arithmetic at the leaf level, yet the hierarchy is as deep as needed to yield the maximal performance gains.

**2.2. Hierarchically semi-separable matrices.** Our primary interest in this section are hierarchically semi-separable (HSS) matrices. These matrices are partitioned in the same way as HODLR matrices but add some additional structure. Given some HODLR partition, with sibling nodes $i$ and $j$ at level $l$, the low-rank block $A(I_i^l, J_j^l)$ can be factored as

$$A(I_i^l, J_j^l) = \tilde{U}_i^{(l)} \tilde{B}_{i,j}^{(l)} \big(\tilde{V}_j^{(l)}\big)^*, \tag{2.2}$$

where $\tilde{U}_i^{(l)} \in \mathbb{R}^{\# I_i^l \times k}$, $\tilde{V}_i^{(l)} \in \mathbb{R}^{\# J_j^l \times k}$ and $\tilde{B}_{i,j}^{(l)} \in \mathbb{R}^{k \times k}$. We call $\tilde{U}_i^{(l)}$ and $\tilde{V}_j^{(l)}$, the generators of $A(I_i^l, J_j^l)$. The core idea of HSS matrices is that they possess an additional hierarchy in their generators. Let $c_1, c_2$ be the children of $i$ and $d_1, d_2$ the children of $j$. For HSS matrices, we can then construct $\tilde{U}_i^{(l)}$ and $\tilde{V}_j^{(l)}$ recursively from the generators of its children:

$$\tilde{U}_i^{(l)} = \begin{bmatrix} \tilde{U}_{c_1}^{(l+1)} & 0 \\ 0 & \tilde{U}_{c_2}^{(l+1)} \end{bmatrix} U_i^{(l)}, \tag{2.3a}$$

$$\tilde{V}_j^{(l)} = \begin{bmatrix} \tilde{V}_{d_1}^{(l+1)} & 0 \\ 0 & \tilde{V}_{d_2}^{(l+1)} \end{bmatrix} V_j^{(l)}, \tag{2.3b}$$

where $U_i^{(l)} \in \mathbb{R}^{2k \times k}$ and $V_j^{(l)} \in \mathbb{R}^{2k \times k}$. In our notation, generators with a $\sim$ denote 'big' matrices with dimensions corresponding to the blocksize of the low-rank block (2.2). To construct them, we apply the nestedness property (2.3) recursively until the leaf level is reached. At the leaf level, we set $U_i^{(L)} = \tilde{U}_i^{(L)}$ and $V_j^{(L)} = \tilde{V}_j^{(L)}$.

This allows us to represent the HSS matrix $A$ in a hierarchical manner. We define

$$B_i^{(l)} = \begin{bmatrix} 0 & \tilde{B}_{c_1,c_2}^{(l+1)} \\ \tilde{B}_{c_2,c_1}^{(l+1)} & 0 \end{bmatrix},$$

$D_i^{(L)} = A(I_i^L, J_i^L)$, as well as

$$U^{(l)} = \operatorname{diag}(U_1^{(l)}, U_2^{(l)}, \ldots, U_{2^l}^{(l)}), \tag{2.4a}$$

$$V^{(l)} = \operatorname{diag}(V_1^{(l)}, V_2^{(l)}, \ldots, V_{2^l}^{(l)}), \tag{2.4b}$$

$$B^{(l)} = \operatorname{diag}(B_1^{(l)}, B_2^{(l)}, \ldots, B_{2^l}^{(l)}), \tag{2.4c}$$

$$D^{(L)} = \operatorname{diag}(D_1^{(L)}, D_2^{(L)}, \ldots, D_{2^L}^{(L)}). \tag{2.4d}$$

Then $A$ can be expressed via the recursion

$$A^{(0)} = B^{(0)}, \tag{2.5a}$$

$$A^{(l)} = U^{(l)} A^{(l-1)} \big(V^{(l)}\big)^* + B^{(l)} \quad \text{for } l = 1, 2, \ldots, L-1 \tag{2.5b}$$

$$A^{(L)} = U^{(L)} A^{(L-1)} \big(V^{(L)}\big)^* + D^{(L)}. \tag{2.5c}$$

We can write out this recursive definition as a telescoping factorization. For a 3-level HSS matrix with balanced cluster trees, this yields

$$A = U^{(3)} \Big( U^{(2)} \Big( U^{(1)} B^{(0)} \big(V^{(1)}\big)^* + B^{(1)} \Big) \big(V^{(2)}\big)^* + B^{(2)} \Big) \big(V^{(3)}\big)^* + \quad D^{(3)} \quad , \tag{2.6}$$

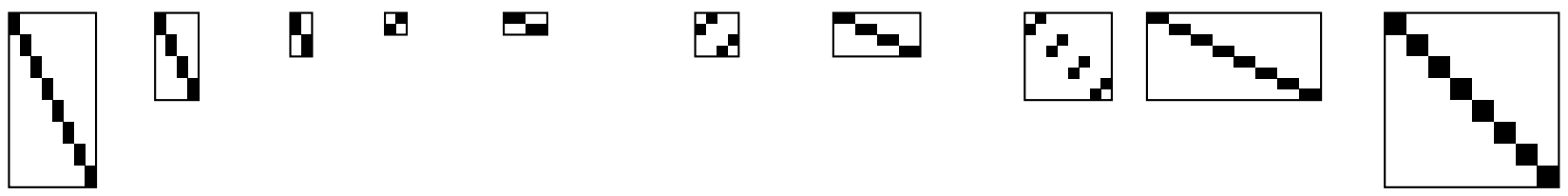

where the structure of each matrix is illustrated below its symbol. A matrix $A$ in HSS format is therefore fully defined by specifying its row- and column-cluster trees $\mathcal{T}_r, \mathcal{T}_c$, as well as the matrices in (2.4).

This telescoping factorization gives rise to an alternative definition of HSS matrices, that is not constructive. We observe in Equation (2.6), that the first row of $A - D^{(3)}$ is spanned by $U_1^{(3)}$. Moreover, we see that this holds for all levels $A^{(l)}$ of the hierarchical definition (2.5), if the diagonal part is disregarded. Consequently, HSS matrices can be defined in the following manner:

Definition 2.3 (HSS matrix). *Consider a matrix $A \in \mathbb{R}^{m \times n}$ with row indices $I$, column indices $J$ and corresponding row- and block-cluster trees $\mathcal{T}_r$, $\mathcal{T}_c$, with identical tree structures.*

- *We select a row partition $I_i^l \in \mathcal{T}_r$ at level $l$. Then the block row $A(I_i^l, J \setminus J_i^l)$, which omits the diagonal part is called a HSS block row. Similarly, we call $A(I \setminus I_\mu^l, J_\mu^l)$ a HSS block column.*
- *We call $A$ an HSS matrix with respect to $\mathcal{T}_r$, $\mathcal{T}_c$, if there exists a small, positive integer $k$, such that the rank of every HSS block row and block column is smaller than, or equal to $k$:*

$$\forall I_i^l \in \mathcal{T}_r : \quad \operatorname{rank} A(I_i^l, J \setminus J_i^l) \le k,$$
$$\forall J_j^l \in \mathcal{T}_c : \quad \operatorname{rank} A(I \setminus I_j^l, J_j^l) \le k.$$

**2.3. Adaptive compression via random sampling.** One of the main limitations of HSS matrices arises from the cost of compressing a general matrix $A \in \mathbb{R}^{n \times n}$ into HSS format. While HSS matrix arithmetic is mostly linear in complexity, a general purpose algorithm for compression relies on the SVD and therefore has quadratic complexity [43]. One remedy is found in randomized algorithms [24], that exploit the fact that the matrix-vector product $x \to Ax$ oftentimes can be computed in linear time [34, 31].

Throughout this paper, we rely on a variant of Martinsson's compression algorithm [34], which is presented for a symmetric matrix $A$ in Appendix A. The original algorithm is modified to sample $p$ extra matrix-vector products of $A$, which are then used to compute the estimator

$$\frac{1}{r}\|(A - A_{\text{HSS}})\tilde{\Omega}\|_{\text{F}}^2, \tag{2.7}$$

where $\tilde{\Omega} \in \mathbb{R}^{n\times r}$ is a Gaussian random matrix. This estimator is used to estimate and control the compression error, such that

$$\|A - A_{\text{HSS}}\|_{\text{F}} \le \epsilon_{\text{HSS}}$$

is met with high probability. Indeed, one can show that the expected value of (2.7) is the norm of $A - A_{\text{HSS}}$ error [24]. If the estimator is larger than the tolerance $\epsilon_{\text{HSS}}$, we can reuse the $r$ sampled matrix-vector products and increase the rank estimate to $k + r$. Martinsson's compression algorithm is then rerun with the new rank estimate.

The advantage of Algorithm A.1 is its $\mathcal{O}(k^2 n)$ complexity, achieved if the matrix-vector products $x \to Ax$, $x \to A^* x$ can be computed in $\mathcal{O}(n)$ time and entries $A(i, j)$ can be accessed in constant time $\mathcal{O}(1)$. This is a considerable improvement over the quadratic complexity of the general purpose algorithm [43], which is based on computing the SVD of individual HSS blocks. It is worth noting however, that Appendix A can often produce a higher approximation error with the same HSS rank. This is why it is essential to use an estimator such as (2.7), to control the error of the compression.

**3. Constructing the preconditioner.** Given a sparse stiffness matrix $A \in \mathbb{R}^{n\times n}$, derived from a Galerkin-type discretization of a PDE, we demonstrate how an approximate factorization

$$A \approx P = LDR, \tag{3.1}$$

can be formed efficiently in $\mathcal{O}(n \log n)$ time. Here, $L$ and $R$ are lower- and upper-block-triangular and $D$ is block-diagonal. It is of equal importance that $P^{-1} = R^{-1}D^{-1}L^{-1}$ can be applied efficiently in $\mathcal{O}(n)$ time. This proposition is based on the hierarchical elimination of degrees of freedom using HSS arithmetic and the proposition that Schur complements can be well-approximated using HSS-matrices. The latter has been extensively discussed [3] and exploited [33, 44, 20, 19] in the existing literature. In this section, we focus on the algorithmic aspect of forming the approximate factorization.

**3.1.** ***LDR*** **factorization.** Let us discuss how an exact factorization $A = LDR$ can be computed using Gaussian elimination and a nested dissection of the computational domain $\Omega$. Figure 1 depicts a nested dissection on the left-hand side. Its hierarchy is represented by its elimination tree $\mathcal{E}$, which is depicted on the right. Each node $\mu$ is associated with a number of degrees of freedom $I^\mu$ and therefore also called a supernode. The elimination tree $\mathcal{E}$ guides the structured elimination of the degrees of freedom of $A$. A node can only be eliminated once all its children nodes have been eliminated. Definition 3.1 summarizes the properties of a valid elimination tree.

Definition 3.1 (Elimination tree). *Let $I = \{1, 2, \ldots, n\}$ be the index set of $A$ and let $\mathcal{E}$ be a tree. Furthermore, let all nodes $\mu \in \mathcal{E}$ be associated to a set of degrees of freedom $I^\mu \subset I$. $\mathcal{E}$ is called an elimination tree of $A$, if*

- *all index sets are disjoint: $\forall \mu \neq \nu : I^\mu \cap I^\nu = \emptyset$*
- *the union of all index sets is the entire set, $\bigcup_{\mu \in \mathcal{E}} I^\mu = I$*
- *For any two indices $i \in I^\mu, j \in I^\nu$ with $\mu > \nu$, $A(i,j) \neq 0$ or $A(j,i) \neq 0$ implies that $\nu$ is a descendant of $\mu$. In other words, the node $\nu$ is included in the subtree rooted at node $\mu$.*

The last property of Definition 3.1 guarantees that all degrees of freedom contributing to a fill-in in $I^\mu$ have already been summed, when $I^\mu$ is being eliminiated. We choose to enumerate the nodes in $\mathcal{E}$ in a post-order manner, as illustrated in Figure 1. To compute the factorization, nodes are eliminated from the bottom-up, following the post-order enumeration. It is worth noting that other ways of traversing $\mathcal{E}$ are also viable. This can be especially useful if we intend to parallelize the algorithm.

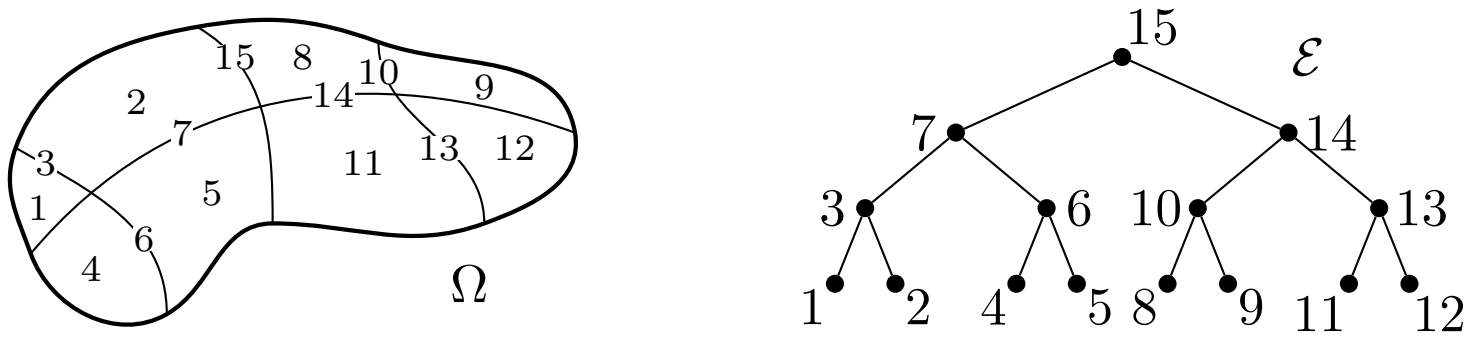

Fig. 1. *Illustration of a nested dissection of the computational domain $\Omega$ using separators. The figure on the right depicts the corresponding post-ordered tree datastructure $\mathcal{E}$, which induces an elimination order.*

Let us assume that we are eliminating the node $\sigma$. We introduce the boundary degrees of freedom

$$B^\sigma = \Big\{ i \in \bigcup\nolimits_{\iota \in ancestors(\sigma)} I^\iota \mid \exists j \in I^\sigma \text{ such that } A(i,j) \neq 0 \text{ or } A(j,i) \neq 0 \Big\},$$

which contains all degrees of freedom in the ancestors that interact with the degrees of freedom in $I^\sigma$. We use the short-hand notation $\hat{A}^{(\sigma)}_{ib} = A(I^\sigma, B^\sigma)$ using subscripts '$i$' and '$b$' to address interior or boundary degrees of freedom of the corresponding node. If $\sigma$ is a leaf node, it suffices to extract the submatrix

$$\hat{A}^{(\sigma)} = \begin{bmatrix} \hat{A}^{(\sigma)}_{ii} & \hat{A}^{(\sigma)}_{ib} \\ \hat{A}^{(\sigma)}_{bi} & \hat{A}^{(\sigma)}_{bb} \end{bmatrix} = \begin{bmatrix} I & 0 \\ -L^{(\sigma)} & I \end{bmatrix} \begin{bmatrix} \hat{A}^{(\sigma)}_{ii} & 0 \\ 0 & S^{(\sigma)} \end{bmatrix} \begin{bmatrix} I & -R^{(\sigma)} \\ 0 & I \end{bmatrix} = \hat{L}^{(\sigma)} \begin{bmatrix} \hat{A}^{(\sigma)}_{ii} & 0 \\ 0 & S^{(\sigma)} \end{bmatrix} \hat{R}^{(\sigma)} \tag{3.2}$$

and factor it. Here, $\hat{L}^{(\sigma)}$ and $\hat{R}^{(\sigma)}$ are block lower and block upper triangular matrices. In practice, it suffices to store the blocks

$$L^{(\sigma)} = -\hat{A}^{(\sigma)}_{bi} \big(\hat{A}^{(\sigma)}_{ii}\big)^{-1}, \tag{3.3a}$$

$$R^{(\sigma)} = -\big(\hat{A}^{(\sigma)}_{ii}\big)^{-1} \hat{A}^{(\sigma)}_{ib}, \tag{3.3b}$$

as we are only interested in applying the action of the inverse of the factors $b \to \big(\hat{L}^{(\sigma)}\big)^{-1} b$ and $b \to \big(\hat{R}^{(\sigma)}\big)^{-1} b$ to a vector $b \in \mathbb{R}^{\#I^\sigma + \#B^\sigma}$. The inverse can easily be applied by changing the sign of the blocks (3.3) and applying the resulting matrices.

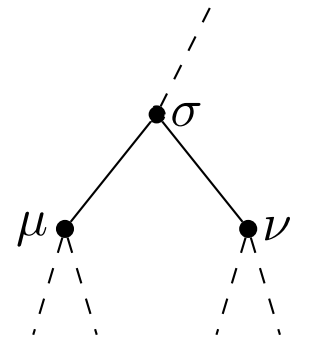

Fig. 2. *Elimination of node $\sigma$. Contributions from the elimination of its children $\mu$ and $\nu$ have to be accounted for.*

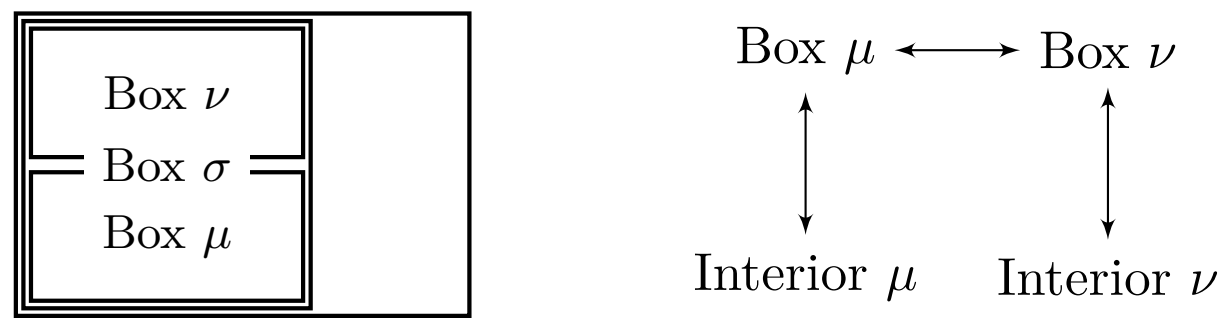

Fig. 3. *Illustration of nested dissection partitioning using boxes. The figure on the right depicts the connectivity of the two children nodes $\mu$ and $\nu$. Interior degrees of freedom are well-separated as they only interact with the boundary of the box.*

The degrees of freedom in $I^\sigma$ have now been decoupled. To continue the factorization, we need to remember contributions to degrees of freedom in $B^\sigma$. The Schur complement

$$S^{(\sigma)} = \hat{A}_{bb}^{(\sigma)} - \hat{A}_{bi}^{(\sigma)} \big(\hat{A}_{ii}^{(\sigma)}\big)^{-1} \hat{A}_{ib}^{(\sigma)}$$

contains all contributions to subsequent nodes and is the only information that has to be passed on.

Let us now consider a non-leaf node as depicted in Figure 2, where all nodes up to the node $\sigma$ including its children $\mu$ and $\nu$ have been factored. So far, the assumptions made were related to the structured Gaussian elimination for sparse matrices. Our next assumption is specific to our algorithm and greatly simplifies its construction. We will assume that the elimination tree $\mathcal{E}$ has the property of *well-separated* nodes. This means that sibling nodes $\mu$ and $\nu$ only interact with each other through disjoint boundaries $B^\mu \cap B^\nu = \emptyset$. In other words, we assume that the elimination of $\mu$ modifies only entries that are left unmodified by the elimination of $\nu$ and vice-versa. Such an elimination tree can be generated by partitioning the grid using boxes as illustrated in Figure 3. The degrees of freedom involved in the elimination of $\sigma$ are $I^\sigma$ and $B^\sigma$. $I^\sigma$ are the interior degrees of freedom marked for elimination at the current step. Then we have the boundary $B^\sigma$, which contains all degrees of freedom that will be modified by the elimination of $\sigma$. Because children nodes are well-separated, $B^\mu$ and $B^\nu$ are disjoint and we can form the four disjoint index sets

$$I^\sigma \cap B^\mu,\ I^\sigma \cap B^\nu,\ B^\sigma \cap B^\mu,\ B^\sigma \cap B^\nu,$$

which form a partitioning of the indices associated with $\sigma$. Using this partitioning, we can now assemble the matrix $\hat{A}^{(\sigma)}$, which contains all the information required for the elimination of $I^\sigma$. Using the partitioning above, which introduces a re-ordering,

we have

$$\hat{A}^{(\sigma)} = \begin{bmatrix} S^{(\mu)}_{ii} & \tilde{A}^{(\mu,\nu)}_{ii} & S^{(\mu)}_{ib} & \tilde{A}^{(\mu,\nu)}_{ib} \\ \tilde{A}^{(\nu,\mu)}_{ii} & S^{(\nu)}_{ii} & \tilde{A}^{(\nu,\mu)}_{ib} & S^{(\nu)}_{ib} \\ S^{(\mu)}_{bi} & \tilde{A}^{(\mu,\nu)}_{bi} & S^{(\mu)}_{bb} & \tilde{A}^{(\mu,\nu)}_{bb} \\ \tilde{A}^{(\nu,\mu)}_{bi} & S^{(\nu)}_{bi} & \tilde{A}^{(\nu,\mu)}_{bb} & S^{(\nu)}_{bb} \end{bmatrix}. \tag{3.4}$$

The Schur complements $S^{(\mu)}$ and $S^{(\nu)}$ were previously formed when factoring the children nodes $\mu$ and $\nu$. With a slight abuse of notation, they are addressed with respect to the indices of their parent node $\sigma$. Thus, $S^{(\mu)}_{ib}$ indicates $S^{(\mu)}(I^\sigma \cap B^\mu, B^\sigma \cap B^\mu)$ and so on. Similarly, we write $\tilde{A}^{(\mu,\nu)}_{ib} = A(I^\sigma \cap B^\mu, B^\sigma \cap B^\nu)$ to indicate entries of the original matrix. With this partitioning in place, we identify the same block structure as in Equation (3.2). As in (3.2), factorization of $\sigma$ requires us to form the Schur complement

(3.5)

$$\begin{aligned} S^{(\sigma)} &= \begin{bmatrix} S^{(\mu)}_{bb} & \tilde{A}^{(\mu,\nu)}_{bb} \\ \tilde{A}^{(\nu,\mu)}_{bb} & S^{(\nu)}_{bb} \end{bmatrix} - \begin{bmatrix} S^{(\mu)}_{bi} & \tilde{A}^{(\mu,\nu)}_{bi} \\ \tilde{A}^{(\nu,\mu)}_{bi} & S^{(\nu)}_{bi} \end{bmatrix} \begin{bmatrix} S^{(\mu)}_{ii} & \tilde{A}^{(\mu,\nu)}_{ii} \\ \tilde{A}^{(\nu,\mu)}_{ii} & S^{(\nu)}_{ii} \end{bmatrix}^{-1} \begin{bmatrix} S^{(\mu)}_{ib} & \tilde{A}^{(\mu,\nu)}_{ib} \\ \tilde{A}^{(\nu,\mu)}_{ib} & S^{(\nu)}_{ib} \end{bmatrix} \\ &= \hat{A}^{(\sigma)}_{bb} - \hat{A}^{(\sigma)}_{bi} \big(\hat{A}^{(\sigma)}_{ii}\big)^{-1} \hat{A}^{(\sigma)}_{ib} \end{aligned}$$

to pass on the information to the next level in the elimination tree. Moreover, we also need to compute and store the inverse of the block matrix $\hat{A}^{(\sigma)}_{ii}$, as well as the left and right factors (3.3). Finally, if $\sigma$ is a the root node, the matrix (3.4) reduces to $\hat{A}^{(\sigma)} = \hat{A}^{(\sigma)}_{ii}$. As such, it is already decoupled from the remaining degrees of freedom. As such, it suffices to compute its inverse, which completes the factorization.

**Algorithm 3.1** Apply the inverse $A^{-1} = R^{-1}D^{-1}L^{-1}$ to a vector $b$

**for all** nodes $\sigma \in \mathcal{E}$ from the bottom up **do**
  Apply $b(B^\sigma) \leftarrow L^{(\sigma)} b(I^\sigma) + b(B^\sigma)$
**end for**
**for all** nodes $\sigma \in \mathcal{E}$ **do**
  Apply $b(I^\sigma) \leftarrow \big(\hat{A}^{(\sigma)}_{ii}\big)^{-1} b(I^\sigma)$
**end for**
**for all** nodes $\sigma \in \mathcal{E}$ from the top down **do**
  Apply $b(I^\sigma) \leftarrow b(I^\sigma) + R^{(\sigma)} b(B^\sigma)$
**end for**
**return** $b$

With the elimination of all nodes, we have implicitly formed the factorization $A = LDR$. We are now interested in applying the inverse $A^{-1} = R^{-1}D^{-1}L^{-1}$. The action of $b \to L^{-1}b$ can be computed by traversing the elimination tree $\mathcal{E}$ from the bottom up and applying the left factor $\hat{L}^{(\sigma)}$ for each node $\sigma$ to the corresponding degrees of freedom. In a similar fashion, we can compute the action of $b \to R^{-1}b$ by traversing $\mathcal{E}$ from the top to the bottom and by applying $\hat{R}^{(\sigma)}$ at each step. The inverse $D^{-1}$ is applied by applying the inverse $\big(\hat{A}^{(\sigma)}_{ii}\big)^{-1}$ locally at each node. This can be done fully in parallel as all nodes are decoupled at this point. The procedure for applying the inverse is summarized in Algorithm 3.1.

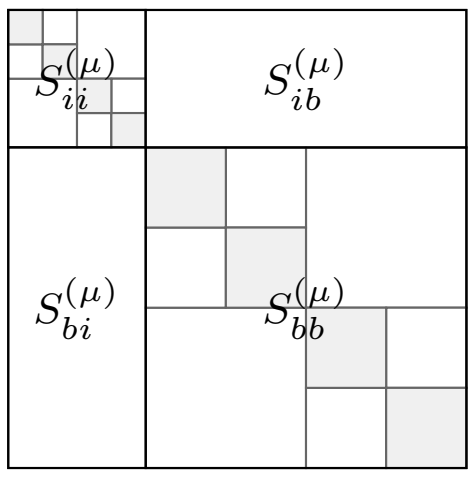

FIG. 4. *Top-level HSS block structure of $S^{(\mu)}$. The partitioning exposes the block $S_{ii}^{(\mu)}$ as a HSS matrix. This block corresponds to degrees of freedom on the interior of the parent of $\mu$, which will be eliminated in the next step of the factorization.*

**3.2. Compression of Schur complements.** So far, all arithmetic has been carried out exactly and we have thus outlined the construction of a sparse direct solver [10]. It is evident that the cost of forming the factorization is driven by the cost of forming Schur complements (3.5) and the inversion of the blocks $\hat{A}_{ii}^{(\sigma)}$, both of which entail dense matrix arithmetic.

In the following sections we will discuss how HSS arithmetic can be leveraged to compute the approximate factorization (3.1) following the hierarchical factorization outlined in the previous section. To make the algorithm matrix-free and efficient, we will switch to HSS arithmetic once a certain switching level $L_{\text{HSS}}$ is reached in $\mathcal{E}$. From this level on, the Schur complements (3.5) are compressed and stored as HSS matrices. This is accomplished via Algorithm A.1, which requires the efficient computation of matrix-vector products $x \to S^{(\sigma)}x$ and $x \to \left(S^{(\sigma)}\right)^* x$, as well as access to individual entries $S^{(\sigma)}(i,j)$, without forming the Schur complement (3.5) explicitly.

Before we continue, however, let us discuss how Schur complements are stored such that the submatrices arising in Equation (3.5) are readily available. The computation of $S^{(\sigma)}$ requires the extraction of the submatrices $S_{ii}^{(\cdot)}$, $S_{bb}^{(\cdot)}$, $S_{ib}^{(\cdot)}$ and $S_{bi}^{(\cdot)}$ from the already compressed Schur complements of its children $S^{(\mu)}$ and $S^{(\nu)}$.

Accordingly, we choose to store the Schur complement of children nodes $\mu$ using an HSS clustering that conforms to the partitioning $I^\sigma \cap B^\mu$, $B^\sigma \cap B^\mu$ of the parent node $\sigma$. This is illustrated in Figure 4. This exposes the degrees of freedoms that will be eliminated next in the top left block and greatly simplifies the extraction of these blocks. Consequently, we can directly access $S_{ii}^{(\mu)}$ and $S_{bb}^{(\mu)}$ as HSS matrices. At the same time, we can extract $S_{ib}^{(\mu)}$ and $S_{bi}^{(\mu)}$ as low-rank matrices with their respective generators.

**3.3. Block inversion of $\hat{A}_{ii}^{(\sigma)}$.** The formation of the parent Schur complement $S^{(\sigma)}$ and the factors $L^{(\sigma)}$ and $R^{(\sigma)}$ requires the efficient application of the inverse

$$\left(\hat{A}_{ii}^{(\sigma)}\right)^{-1} = \begin{bmatrix} S_{ii}^{(\mu)} & \tilde{A}_{ii}^{(\mu,\nu)} \\ \tilde{A}_{ii}^{(\nu,\mu)} & S_{ii}^{(\nu)} \end{bmatrix}^{-1}. \tag{3.6}$$

to the right $X \to \left(\hat{A}_{ii}^{(\sigma)}\right)^{-1} X$ and to the left $X \to X\left(\hat{A}_{ii}^{(\sigma)}\right)^{-1}$.

Given an invertible block matrix

$$B = \begin{bmatrix} B_{11} & B_{12} \\ B_{21} & B_{22} \end{bmatrix},$$

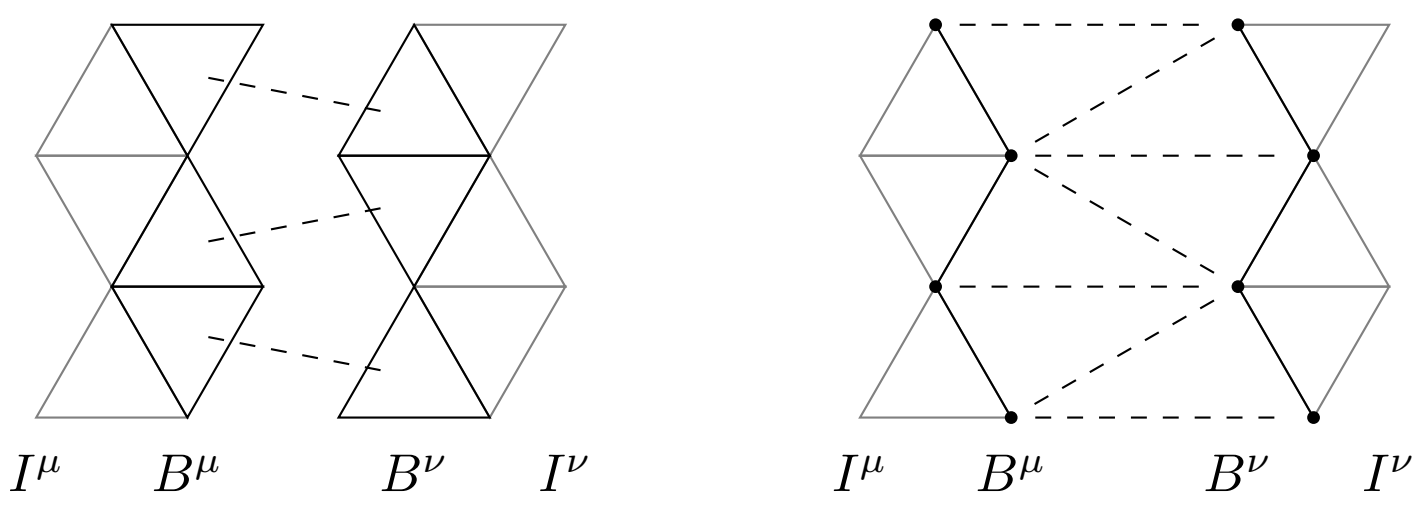

FIG. 5. *Connectivity between the children elements $\mu$ and $\nu$ in the interior of $\sigma$. The figure on the left depicts the situation for discontinuous Galerkin discretizations. The figure on the right depicts the connectivity arising from continuous Galerkin discretizations. Dashed lines represend interactions and entries in $\tilde{A}_{ii}^{(\mu,\nu)}$ and $\tilde{A}_{ii}^{(\nu,\mu)}$.*

a partitioned right-hand side $X = \begin{bmatrix} X_1 \\ X_2 \end{bmatrix}$, with partitioning corresponding to the blocks in $B$, we can apply the inverse to the right by using Algorithm 3.2. A similar algorithm

**Algorithm 3.2** Right-apply block inverse of $B$ to $X$

Partition $X = \begin{bmatrix} X_1 \\ X_2 \end{bmatrix}$, corresponding to the blocks in $B$
$X_1 \leftarrow B_{11}^{-1} X_1$
$X_2 \leftarrow X_2 - B_{21} X_1$
Form the Schur complement $\tilde{S} = B_{22} - B_{21} B_{11}^{-1} B_{12}$
$X_2 \leftarrow \tilde{S}^{-1} X_2$
$X_1 \leftarrow X_1 - B_{11}^{-1} B_{12} X_2$
**return** $X = \begin{bmatrix} X_1 \\ X_2 \end{bmatrix}$

can be constructed to apply the block inverse $B^{-1}$ to the left. It is evident that a efficient computation would require efficient methods to form the Schur complement $S_{22}$ and invert $B_{11}$, $B_{22}$ and $S_{22}$.

We assume that the Schur complements $S^{(\mu)}$, $S^{(\nu)}$ have already been compressed, which makes the blocks $S_{ii}^{(\mu)}$ and $S_{ii}^{(\nu)}$ readily available in HSS form. The off-diagonal blocks $\tilde{A}_{ii}^{(\mu,\nu)}$ and $\tilde{A}_{ii}^{(\nu,\mu)}$ are submatrices of the original matrix and therefore sparse. They represent interactions between the two children nodes $\mu$ and $\nu$ within the interior of $\sigma$. Figure 5 depicts the corresponding connectivity for discontinuous Galerkin and continuous Galerkin discretizations.

Let us briefly discuss why these matrices can be represented by hierarchical matrices with low off-diagonal ranks as this will allow us to construct an accelerated algorithm for the block inverse. To this end, we consider the matrix $\tilde{A}_{ii}^{(\mu,\nu)}$, and the row cluster tree $\mathcal{T}_r^{(\mu)}$ of $S_{ii}^{(\mu)}$ along with the column cluster tree $\mathcal{T}_c^{(\nu)}$ of $S_{ii}^{(\nu)}$. We recall Definition 2.3, which states that if every HSS block row and HSS block column is low-rank, the matrix can be approximated by a HSS matrix.

In most discontinuous Galerkin type discretizations, elements only interact across faces as in Figure 5. As a consequence, we can keep elements intact when the nested dissection is performed. Therefore each element in $B^\mu$ only interacts with one distinct element in $B^\nu$. If we number these elements so they have the same position in their respective cluster trees $\mathcal{T}_r^{(\mu)}$ and $\mathcal{T}_c^{(\nu)}$, we have

$$\forall I_i^l \in \mathcal{T}_r^{(\mu)} : \quad \operatorname{rank} A(I_i^l, J \setminus J_i^l) = 0,$$
$$\forall J_j^l \in \mathcal{T}_c^{(\nu)} : \quad \operatorname{rank} A(I \setminus I_j^l, J_j^l) = 0,$$

which makes $\tilde{A}_{ii}^{(\mu,\nu)}$ an HSS matrix of HSS rank 0, which is coincidentally a block-diagonal matrix. As previously mentioned, this requires some care in the construction of the nested dissection $\mathcal{E}$ and the numbering of degrees of freedom within. Moreover, the blocks in $I_i^L$ at the leaf level should conform to the elements of the discontinuous Galerkin discretization. In practice, this can be done by choosing a HSS block size $\beta$, which corresponds to a multiple of the number of degrees of freedom per element. This also makes the algorithm robust in the case that the numbering is not conforming.

For continuous Galerkin discretizations, the situation is more complex due to degrees of freedom being shared among elements. To ensure that boxes are well-separated, we split elements at the interfaces as seen in Figure 5. Already for first order finite elements, it is not possible to group degrees of freedom of $B^\mu$ and $B^\nu$ into clusters, such that each cluster of $B^\mu$ only interacts with one cluster of $B^\nu$ and vice versa. However, we can form clusters such that each cluster mainly interacts with one cluster and has only one additional connection to another cluster. Consequently, we have

$$\forall I_i^l \in \mathcal{T}_r^{(\mu)} : \quad \operatorname{rank} A(I_i^l, J \setminus J_i^l) \le \operatorname{nnz} A(I_i^l, J \setminus J_i^l) = 1,$$
$$\forall J_j^l \in \mathcal{T}_c^{(\nu)} : \quad \operatorname{rank} A(I \setminus I_j^l, J_j^l) \le \operatorname{nnz} A(I \setminus I_j^l, J_j^l) = 1,$$

which makes $\tilde{A}_{ii}^{(\mu,\nu)}$ an HSS matrix of rank 1. As in the discontinuous Galerkin case, we have to make sure that the nested dissection reordering numbers degrees of freedom accordingly.

Hence, $\tilde{A}_{ii}^{(\mu,\nu)}$ and $\tilde{A}_{ii}^{(\nu,\mu)}$ can indeed be expressed as HSS matrices with low HSS rank. As a consequence, we modify Algorithm 3.2 using HSS arithmetic to form the intermediate Schur complement $\tilde{S}$ in HSS form. This requires the compression of $\tilde{A}_{ii}^{(\mu,\nu)}$ and $\tilde{A}_{ii}^{(\nu,\mu)}$ with compatible row and column clusters among $S_{ii}^{(\nu)}$, $\tilde{A}_{ii}^{(\nu,\mu)}$, $S_{ii}^{(\mu)}$ and $\tilde{A}_{ii}^{(\mu,\nu)}$. The procedure is described in Algorithm 3.3. The compression of $\tilde{A}_{ii}^{(\mu,\nu)}$ and $\tilde{A}_{ii}^{(\nu,\mu)}$ can be realized in linear complexity via Algorithm A.1. In practice, most entries lie in the block diagonal and the compression can be accelerated by computing it directly from the sparse representation. The formation of $A^{-1}B$ in HSS format is based on algorithms presented in [7, 36]. This algorithm has a complexity of $\mathcal{O}(k^2 n)$, where $n$ is the corresponding matrix size and $k$ the maximum HSS rank encountered. Thus the overall complexity of Algorithm 3.3 is equally $\mathcal{O}(k^2 n)$.

**3.4. Computing the factors $L^{(\sigma)}$ and $R^{(\sigma)}$.** We consider the computation of the factors $L^{(\sigma)}$ and $R^{(\sigma)}$. The right factor reads

$$R^{(\sigma)} = -\big(\hat{A}_{ii}^{(\sigma)}\big)^{-1} \hat{A}_{ib}^{(\sigma)} = -\Big(\hat{A}_{ii}^{(\sigma)}\Big)^{-1} \begin{bmatrix} S_{ib}^{(\mu)} & \tilde{A}_{ib}^{(\mu,\nu)} \\ \tilde{A}_{ib}^{(\nu,\mu)} & S_{ib}^{(\nu)} \end{bmatrix}. \tag{3.7}$$

We already know how to apply the inverse of $\tilde{A}_{ii}^{(\sigma)}$ efficiently to a vector or matrix with compatible block structure. A closer look at $\tilde{A}_{ib}^{(\sigma)}$ reveals that the diagonal blocks $S_{ib}^{(\mu)}$ and $S_{ib}^{(\nu)}$ have a low-rank structure, while the off-diagonal blocks $\tilde{A}_{ib}^{(\mu,\nu)}$, $\tilde{A}_{ib}^{(\nu,\mu)}$

**Algorithm 3.3** Right-apply the block inverse of $\hat{A}_{ii}^{(\sigma)}$ using HSS arithmetic

Partition $X = \begin{bmatrix} X_1 \\ X_2 \end{bmatrix}$, corresponding to the blocks in $\hat{A}_{ii}^{(\sigma)}$
Extract row and column clusters $\mathcal{T}_r^{(\mu)}$, $\mathcal{T}_c^{(\mu)}$ of $S_{ii}^{(\mu)}$
Extract row and column clusters $\mathcal{T}_r^{(\nu)}$, $\mathcal{T}_c^{(\nu)}$ of $S_{ii}^{(\nu)}$
Compress $\tilde{A}_{ii}^{(\mu,\nu)}$ in HSS form with row and column clusters $\mathcal{T}_r^{(\mu)}$, $\mathcal{T}_c^{(\nu)}$
Compress $\tilde{A}_{ii}^{(\nu,\mu)}$ in HSS form with row and column clusters $\mathcal{T}_r^{(\nu)}$, $\mathcal{T}_c^{(\mu)}$
$X_1 \leftarrow \big(S_{ii}^{(\mu)}\big)^{-1} X_1$
$X_2 \leftarrow X_2 - \hat{A}_{ii}^{(\mu,\nu)} X_1$
Form the Schur complement $\tilde{S} = S_{ii}^{(\nu)} - \tilde{A}_{ii}^{(\nu,\mu)}\big(S_{ii}^{(\mu)}\big)^{-1}\tilde{A}_{ii}^{(\mu,\nu)}$ in HSS form
Compute $X_2 \leftarrow \tilde{S}^{-1} X_2$
$X_1 \leftarrow X_1 - \big(S_{ii}^{(\mu)}\big)^{-1}\hat{A}_{ii}^{(\mu,\nu)} X_2$
**return** $X = \begin{bmatrix} X_1 \\ X_2 \end{bmatrix}$

are sparse. We have

$$\begin{aligned} \operatorname{rank} R^{(\sigma)} &= \operatorname{rank} \hat{A}_{ib}^{(\sigma)} \\ &\leq \operatorname{rank} S_{ib}^{(\mu)} + \operatorname{rank} S_{ib}^{(\nu)} + \operatorname{rank} \tilde{A}_{ib}^{(\mu,\nu)} + \operatorname{rank} \tilde{A}_{ib}^{(\nu,\mu)} \\ &\leq \operatorname{hssrank} S^{(\mu)} + \operatorname{hssrank} S^{(\nu)} + \operatorname{nnz} \tilde{A}_{ib}^{(\mu,\nu)} + \operatorname{nnz} \tilde{A}_{ib}^{(\nu,\mu)} \end{aligned}$$

which implies that $R^{(\sigma)}$ has a low-rank structure itself if $\tilde{A}_{ib}^{(\mu,\nu)}$ and $\tilde{A}_{ib}^{(\nu,\mu)}$ are sufficiently sparse. $\tilde{A}_{ib}^{(\mu,\nu)}$ reflects interactions between the interior degrees of freedom $I^\sigma \cap B^\mu$, associated with box $\mu$ and the boundary degrees of freedom $B^\sigma \cap B^\nu$ associated with box $\nu$. This interaction is typically extremely sparse and will be zero in most cases. In two dimensions, it corresponds to interactions in one point. In three dimensions, this interaction corresponds to interactions along a line. Consequently, the rank of $L^{(\sigma)}$ and $R^{(\sigma)}$ is bounded by double the maximum HSS rank that is encountered among Schur complements plus a small constant. A similar rank estimate can be derived for the left factor $L^{(\sigma)}$. This low-rank property is well-known for matrices arising from Galerkin discretizations of elliptic problems [38, 45]. Here, we have seen how this property is connected to the approximability of Schur complements using hierarchical matrices, another property that has been proven for such matrices [4, 2, 3].

Consequently, we can store the factors $L^{(\sigma)}$ and $R^{(\sigma)}$ in low-rank format represented by their generators. Moreover, as we can apply both $\big(\hat{A}_{ii}^{(\sigma)}\big)^{-1}$ and $\hat{A}_{ib}^{(\sigma)}$ in linear time, we can use a randomized algorithm to compress and store their low-rank representations efficiently in $\mathcal{O}(k^3 n)$ time [24].

**3.5. Summary of the algorithm.** We have now discussed all steps necessary to compute the HSS representation of the Schur complement (3.5) via random sampling. As mentioned in Subsection 3.2, we require an efficient way of computing $x \to S^{(\sigma)}x$, $x \to \big(S^{(\sigma)}\big)^* x$ and $S^{(\sigma)}(i,j)$. The first matrix in (3.5), $\hat{A}_{bb}^{(\sigma)}$ consists of HSS matrices on the diagonal and sparse matrices on the off-diagonal. The cost of applying $x \to \hat{A}_{bb}^{(\sigma)} x$ is $\mathcal{O}(kn + n_{\mathrm{nz}} n)$, where $n$ is the matrix size, $n_{\mathrm{nz}}$ the maximum number of nonzero entries per line and $k$ is the maximum HSS rank. On the other hand, the cost of accessing individual entries is $\mathcal{O}(\log n_{\text{nz}})$, if it lies in the sparse blocks and $\mathcal{O}(k \log n)$ for the HSS blocks. The second matrix to be considered in (3.5) is

$$\hat{A}_{bi}^{(\sigma)} \big(\hat{A}_{ii}^{(\sigma)}\big)^{-1} \hat{A}_{ib}^{(\sigma)} = \hat{A}_{bi}^{(\sigma)} R^{(\sigma)}.$$

As we have previously established in Subsection 3.4, both $\hat{A}_{bi}^{(\sigma)}$ and $R^{(\sigma)}$ are low-rank matrices. As a consequence, we can compute matrix-vector products in $\mathcal{O}(kn)$ time and access individual entries in $\mathcal{O}(k)$ time.

Therefore, we can compress the Schur complement $S^{(\sigma)}$ with Algorithm A.1 in $\mathcal{O}(k^2 n \log n)$ operations[1]. It is worth noting that both the matrix-vector products $x \to S^{(\sigma)}x$, $x \to \big(S^{(\sigma)}\big)^* x$, as well as access to individual entries can easily be adapted to allow for permutations of the indices. As such, we can compress $S^{(\sigma)}$ with the block structure depicted in Figure 4, which is required for the factorization of the parent node.

With the efficient computation of $L^{(\sigma)}$, $R^{(\sigma)}$ and $S^{(\sigma)}$ in place, we have all pieces needed to form the approximate factorization (3.1), which is summarized in Algorithm 3.4.

**Algorithm 3.4** Compute approximate factorization $P = LDR$

```
Generate a suitable elimination tree ℰ
for all nodes σ ∈ ℰ from the bottom up do
  if level(σ) < L_HSS then
    Implicitly form Â_ii^(σ) by extracting the corresponding blocks
    Factor Â_ii^(σ) so that Algorithm 3.3 can be applied
    Compute L^(σ) and R^(σ) in low-rank format
  else
    Factor Â_ii^(σ) using dense arithmetic
    Form L^(σ) and R^(σ) using dense arithmetic
  end if
  if σ is the root node then
    Do nothing
  else if level(σ) = L_HSS then
    Form S^(σ) densely and compress to HSS form
  else if level(σ) < L_HSS then
    Form S^(σ) directly in HSS form via Algorithm A.1
  else
    Form and store S^(σ) as dense matrix
  end if
end for
```

As with the exact factorization, we require an algorithm for applying the approximate inverse $P^{-1} = R^{-1} D^{-1} L^{-1}$. This can be done by simply using Algorithm 3.1, where we switch to HSS arithmetic whenever a matrix is stored in HSS format.

**3.6. Computational complexity.** It is useful to summarize the overall cost of forming the approximate factorization. As all operations above the switching level

[1] The increased cost of accessing individual entries is not a problem for the compression algorithm as only $\mathcal{O}(kn)$ entries have to be accessed. As a result, the overall complexity for the HSS compression is $\mathcal{O}(k^2 n \log n)$.

$L_{\text{HSS}}$ are decoupled, their asymptotic cost is linear as the hierarchy grows proportionally with the overall size of the matrix $n$. As such, we focus on the nested dissection hierarchy above the switching level $L_{\text{HSS}}$.

At each branch node on level $l$, we have to process matrices of size $n_l$. We assume that HSS ranks and ranks of $L^{(\sigma)}$ and $R^{(\sigma)}$ can be bounded by the HSS rank of the top-level Schur complement $k$. Table 1 provides an overview of the cost associated with each operation that is carried out at every node above the switching level. Most notably, we see that the cost is dominated by the compression of the Schur complement, due to the increased cost of accessing individual entries in $S^{(\sigma)}$. Moreover, assembling the blocks and applying them to a vector has a hidden dependency on the number of non-zeo entries per column $n_{\text{nz}}$, which we do not account for, as it is usually a constant for a given problem. It may play a role, e.g. if $p$-refinement is considered. The overall complexity can be determined by summing over all leaves. In

Table 1
*Summary of the computational cost of each operation at a single node*

| operation | operation count |
|---|---|
| **factorization** | |
| form $\hat{A}_{bb}^{(\sigma)}$, $\hat{A}_{bi}^{(\sigma)}$, $\hat{A}_{ii}^{(\sigma)}$ and $\hat{A}_{ib}^{(\sigma)}$ | $\mathcal{O}(kn_l \log n_l)$ |
| factor $\hat{A}_{ii}^{(\sigma)}$ | $\mathcal{O}(k^2 n_l)$ |
| compute $L^{(\sigma)}$ and $R^{(\sigma)}$ | $\mathcal{O}(k^3 n_l)$ |
| compress $S^{(\sigma)}$ | $\mathcal{O}(k^2 n_l \log n_l)$ |
| **application** | |
| apply $\big(\hat{A}_{ii}^{(\sigma)}\big)^{-1}$ | $\mathcal{O}(kn_l)$ |
| apply $L^{(\sigma)}$ and $R^{(\sigma)}$ | $\mathcal{O}(kn_l)$ |
| apply $S^{(\sigma)}$ | $\mathcal{O}(kn_l)$ |

relation to the overall size $n$, the size of matrices to factor at each level grows as

$$n_l \sim 2^{-\frac{d-1}{d}(l-1)} n^{\frac{d-1}{d}}, \tag{3.8}$$

where $d$ is the dimension of the problem [35]. In other words, the size of the separators in the nested dissection is divided by $2^{d-1}$ every $d$ levels. As the number of operations at each node is $\mathcal{O}(k^2 n_l \log n_l + k^3 n_l)$, the overall number of operations is

$$\begin{aligned} W &\sim \sum_{l=1}^{L_{\text{HSS}}} 2^{l-1}(k^3 n_l + k^2 n_l \log n_l) \\ &\lesssim n^{\frac{d-1}{d}}(k^3 + k^2 \log n) \sum_{l=1}^{L_{\text{HSS}}} 2^{\frac{1}{d}(l-1)} \sim n^{\frac{d-1}{d}}(k^3 + k^2 \log n) 2^{\frac{1}{d} L_{\text{HSS}}} \\ &\sim k^2 n \log n + k^3 n, \end{aligned} \tag{3.9}$$

where we have used that $L_{\text{HSS}} \sim \log n$. We conclude that the approximate factorization can be formed in quasilinear time $\mathcal{O}(k^2 n \log n + k^3 n)$, regardless of the dimension $d$, assuming that the ranks $k$ stay constant. A similar argument can be made for the cost of applying the preconditioner, which suggests linear complexity of $\mathcal{O}(kn)$.

## 4. Numerical results.

We present numerical results and demonstrate the viability of our preconditioner. Throughout most experiments, we use a tolerance of $\epsilon_{\text{HSS}} = 10^{-6}$ and a block size of $\beta = 10(p+1)(p+2)/2$ for both HSS compression and arithmetic. This block size corresponds to the number of degrees of freedom in 10 elements, if a discontinuous Galerkin discretization of polynomial degree $p$ is used for a two-dimensional system of PDEs. These choices are based on our own experience and not necessarily optimal. In optimized implementations, the block size should be adapted such that blocks are only compressed if there is a performance benefit.

Furthermore, an adaptive hierarchy is utilized, which becomes deeper under $h$-refinement. In particular, the nested dissection proceeds until there are less than 10 elements left in each box. Unless stated otherwise, we switch to HSS arithmetic 4 levels above the lowest level in the nested dissection hierarchy.

We test our preconditioner using restarted GMRES, where the Arnoldi basis is recomputed every 10 iterations. The performance is measured by considering the relative residual $\left\|P^{-1}Ax_i - P^{-1}b\right\| / \left\|P^{-1}b\right\|$, where $P^{-1}$ represents the action of our preconditioner. This process is terminated once a relative residual smaller than $10^{-9}$ is achieved. Alternatively, we terminate the computation at a maximum of 30 GMRES iterations.

**4.1. Poisson problem.** We consider the Poisson problem in two dimensions. Let $\Omega \subset \mathbb{R}^2$ be our domain of interest with Dirichlet and Neumann boundaries $\Gamma_D$ and $\Gamma_N$, such that $\Gamma_D \cup \Gamma_N = \partial\Omega$ and $\Gamma_D \cap \Gamma_N = \emptyset$. The unit normal vector on the boundary $\partial\Omega$ is denoted $\hat{n}$. The Poisson problem is: Given a right-hand side $f : \Omega \to \mathbb{R}$, find a solution $u : \Omega \to \mathbb{R}$, that satisfies

$$-\nabla^2 u = f \text{ in } \Omega \tag{4.1a}$$

$$u = g_D \text{ on } \Gamma_D \tag{4.1b}$$

$$\hat{n} \cdot \nabla u = g_N \text{ on } \Gamma_N. \tag{4.1c}$$

We consider continuous Galerkin (CG), as well as interior penalty discontinuous Galerkin (IPDG) discretizations of (4.1) on the unit square $\Omega = [-1, 1]^2$ as discussed in [25]. We choose Dirichet boundary conditions $g_D = 0$ on the entire boundary and a constant right-hand side $f = 1$. The domain is discretized with a regular, triangular mesh and the nested dissection is generated by subsequent subdivisions of the mesh. To keep the aspect ratio of individual boxes roughly balanced, we alternate the orientation of the subdivision. Once there are less than 25 elements in each remainder mesh, we terminate the process. The preconditioner is then generated for the resulting stiffness matrices, using the nested dissection elimination tree. Figure 6 shows the relative residual at each GMRES iteration $i$. We can see that the preconditioner works equally well for IPDG and CG discretizations of the Poisson problem. Moreover, the number of GMRES iterations increases only slightly as $1/h$ increases.

**4.2. Helmholtz problem.** Next, we consider IPDG discretizations of the inhomogeneous Helmholtz problem

$$-\nabla^2 u - \kappa^2 u = f \text{ in } \Omega \tag{4.2a}$$

$$u = g_D \text{ on } \Gamma_D \tag{4.2b}$$

$$\nabla u \cdot \hat{n} = g_N \text{ on } \Gamma_N. \tag{4.2c}$$

Here, $\kappa$ is called the wavenumber. If the right-hand side is set to zero, the formulation becomes equivalent to the eigenvalue problem of the Poisson problem, and $\kappa^2$ corresponds to the eigenvalue of the Laplacian. Therefore, if a wavenumber that coincides with the square-root of an eigenvalue of the Laplacian is chosen, the problem becomes singular. It is useful to refer to the wavelength of the problem, which is $\lambda = 2\pi/\kappa$.

For our numerical tests, we keep the square domain with Dirichlet boundary conditions $g_D = 0$ on the entire boundary. The mesh generation remains unchanged

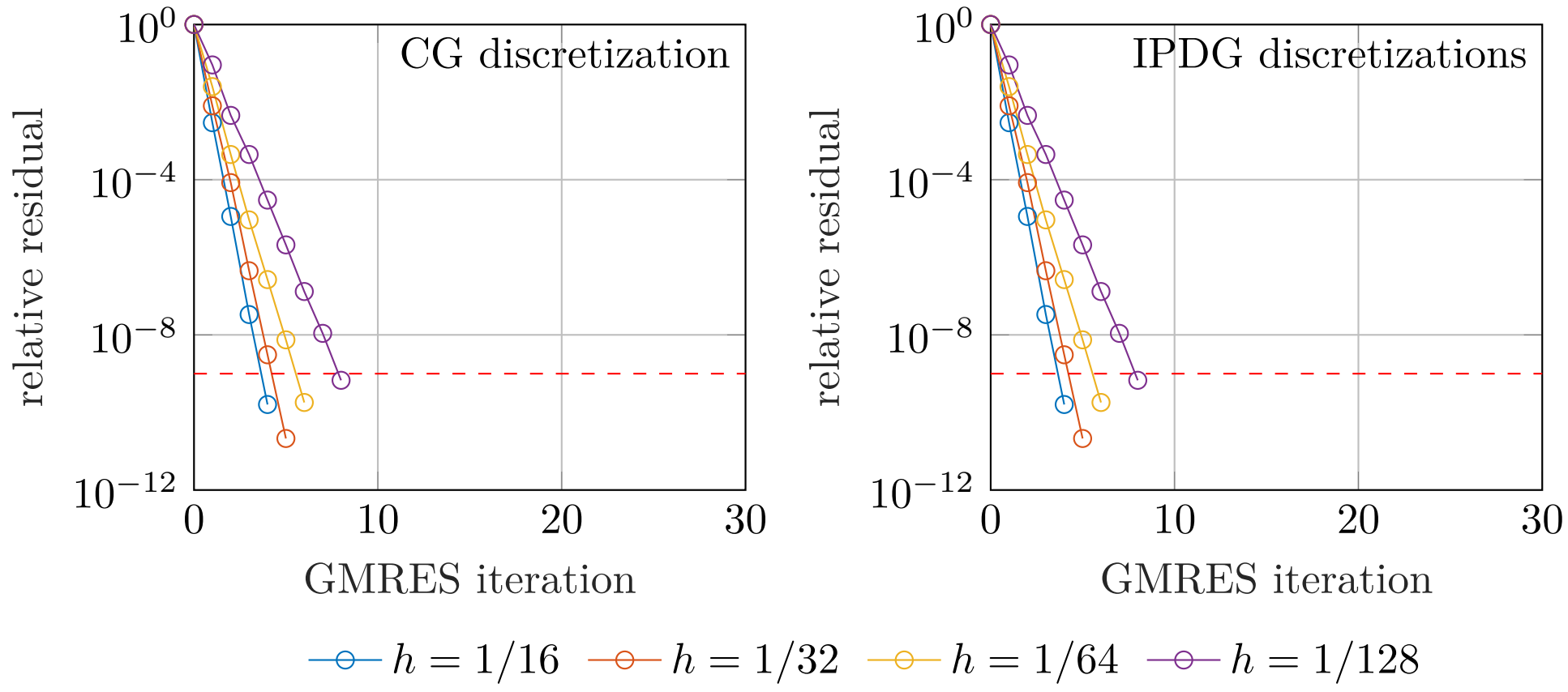

FIG. 6. *Preconditioner performance for the Poisson problem. On the left we depict the performance in the case of $p = 4$ CG discretization. On the right we show the performance of our preconditioner for a $p = 4$ IPDG discretization of the same problem.*

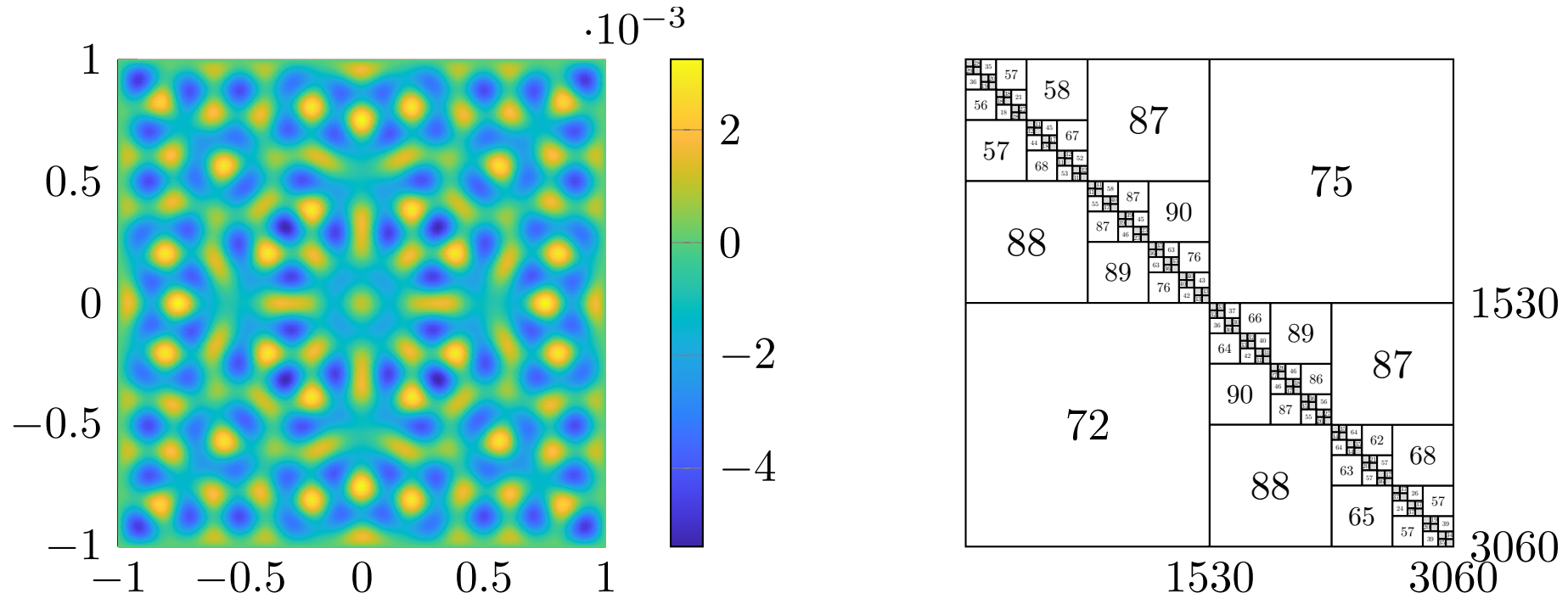

FIG. 7. *On the left we show a solution to (4.2) for $\kappa = 32$, computed with a discontinuous Galerkin appoimation with $h = 1/64$ and $p = 2$. The HSS structure and off-diagonal ranks of the top-most Schur complement are displayed on the right. The smallest off-diagonal blocks have ranks of approximately 35.*

and the nested dissection is performed until a maximum of 10 elements is contained in every leaf level box. The rule of thumb for standard techniques is to have a minimum of 10 grid points per wavelength [29]. Our experiments are well within this limit, even for the wavenumber $\kappa = 64$. Figure 7 depicts the IPDG solution for a polynomial degree of 2, a mesh width of $h = 1/64$ and a wave number of $\kappa = 32$. Again, 4 levels above the leaf level, the preconditioner switches to HSS compression, which amounts to a 7 level hierarchy. On the right Figure 7 we depict the HSS ranks of the top-level Schur complement. We observe that off-diagonal ranks are bounded by 90 and that the Schur complement is well-represented by the HSS format. At the leaf level, the ranks are approximately 30, which is not optimal as the leaf-level boxes have a size of $60 \times 60$. This illustrates that if performance is critical, the HSS leaf level size has to be chosen more carefully to ensure that HSS arithmetic speeds up computations. Figure 8 depicts the residual history for the same discretization with various wavenumbers ranging from $\kappa = 4$ to $\kappa = 64$. For

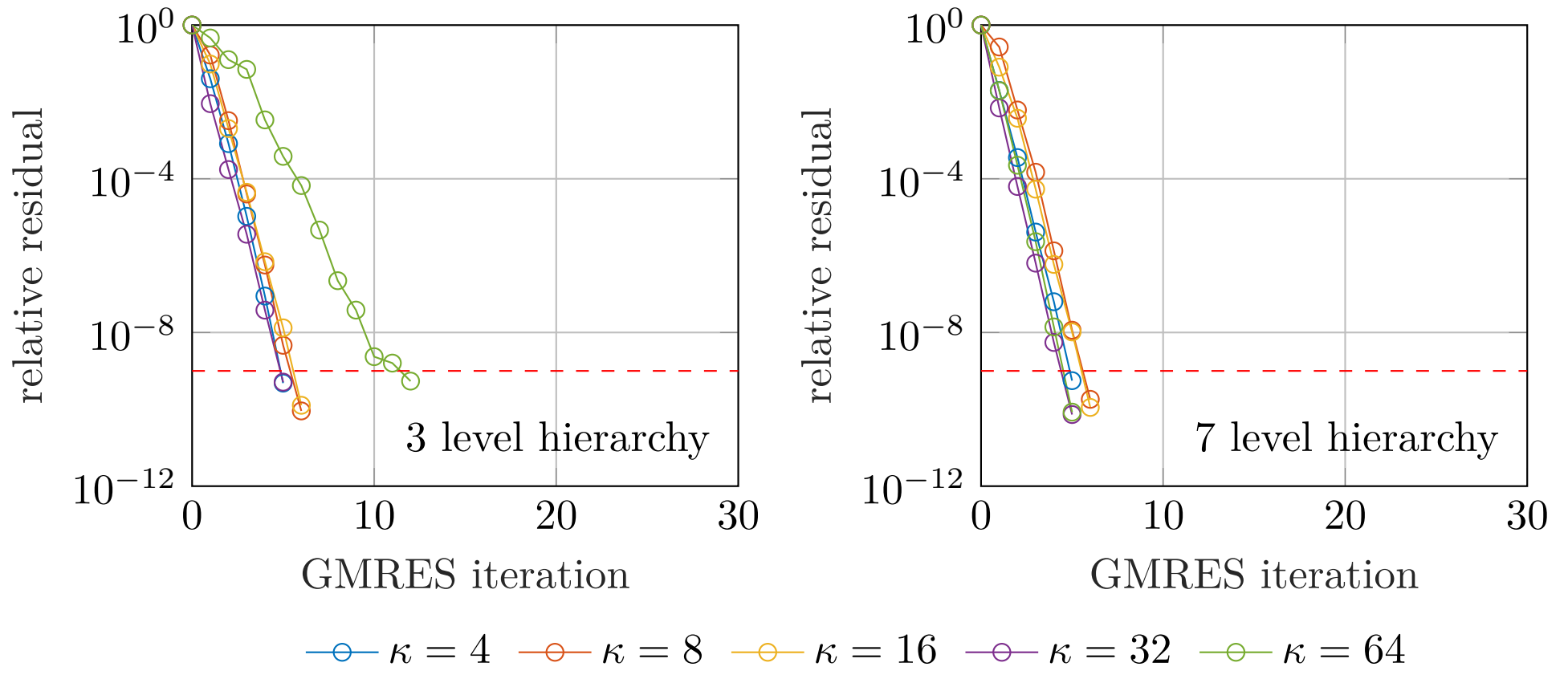

FIG. 8. *Relative residual at each iteration of preconditioned GMRES for the IPDG discretization of the Helmholtz problem (4.2) with $h = 1/64$ and $p = 2$. The test were performed using 3- and 7-level hierarchies and various wave numbers.*

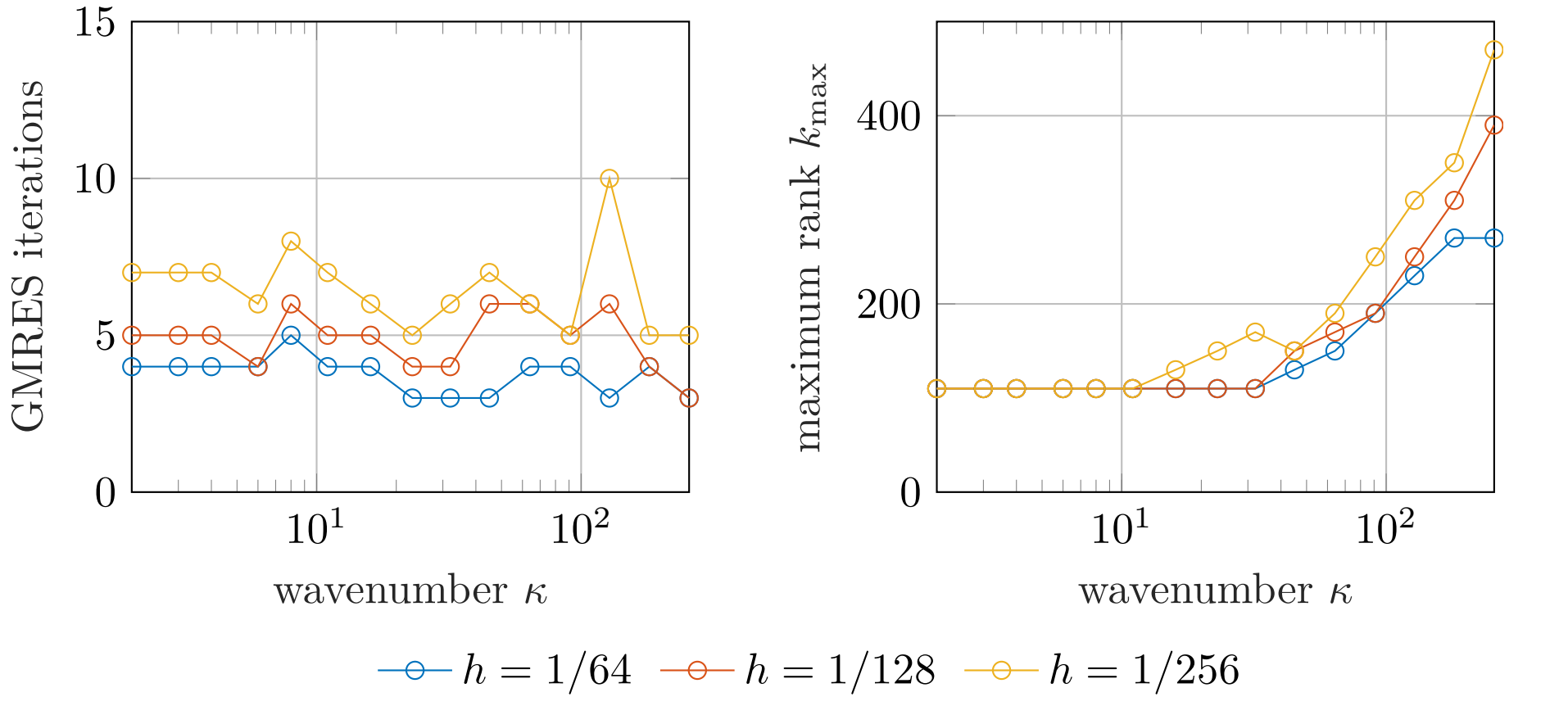

FIG. 9. *Preconditioner performance at exponentially increasing wave numbers and $p = 2$. On the left the number of GMRES iterations to reach a tolerance of $10^{-9}$ are depicted. On the right the maximum among HSS ranks of Schur complements and ranks of Gauss transforms $k_{max}$ is shown.*

all wavenumbers, the preconditioner is able to ensure fast convergence. This includes high wavenumbers, where many other preconditioning techniques typically fail. For a 3 level hierarchy, we can observe that the high wavenumber $\kappa = 64$ leads to a higher number of iterations, whereas in all other cases, the performance is consistent. Figure 9 depicts both the number of GMRES iterations needed to achieve a tolerance of $10^{-9}$ and the maximal rank $k_{\max}$, for exponentially increasing wavenumbers ranging from $\kappa = 2$ to $\kappa = 256$. The maximum rank $k_{\max}$ is the maximum rank encountered across all low-rank matrices of the factorization. This includes both low-rank blocks of Schur complements and the low-rank factors of $L^{(\sigma)}$ and $R^{(\sigma)}$. Typically, it is twice the HSS rank of the top-level Schur complement. Experiments are run for different mesh widths $1/h$ at a polynomial degree of $p = 2$ to ensure that the solution remains

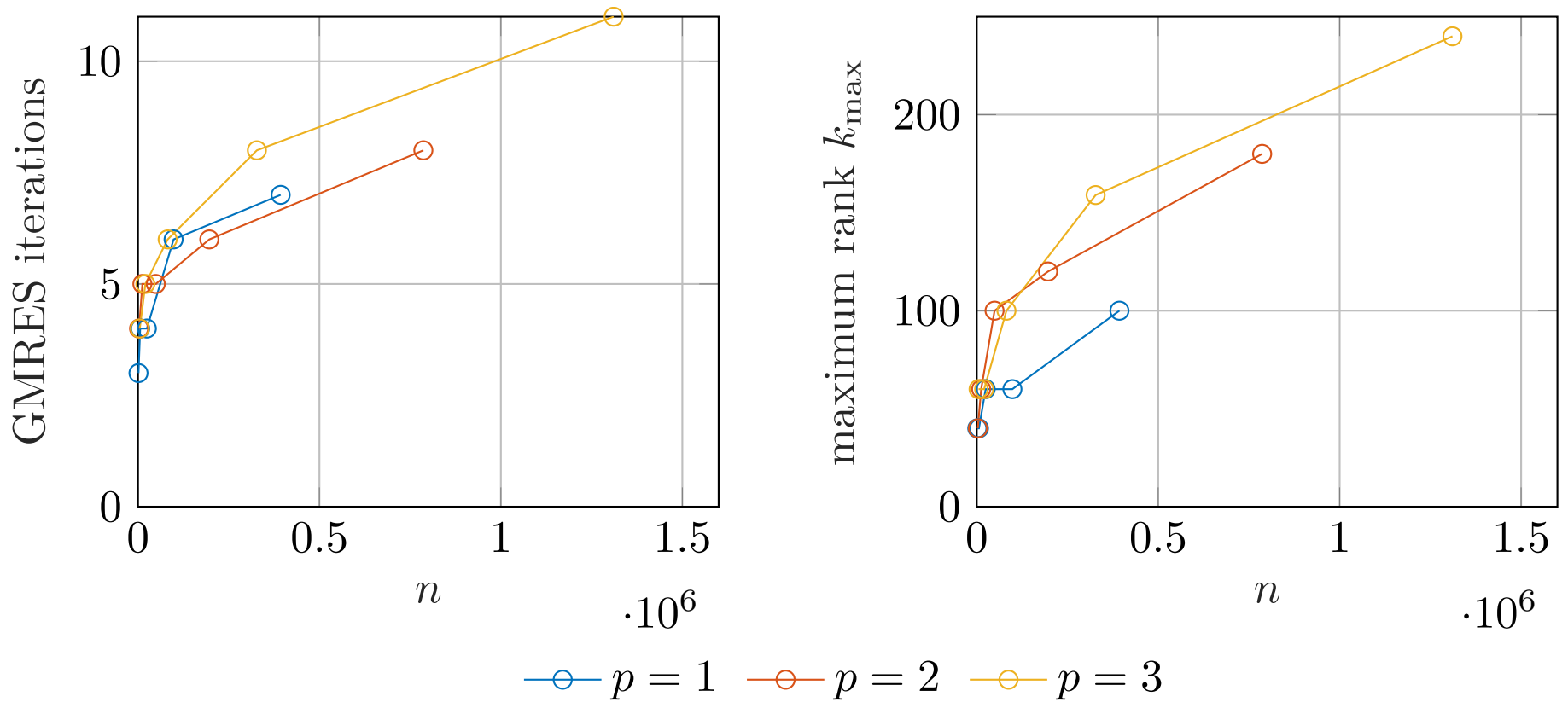

FIG. 10. *Preconditioner performance for the Helmholtz problem under $h$-refinement. $n$ is the number of degrees of freedom.*

well-resolved.

We observe that the amount of iterations required to compute an accurate solution essentially remains constant, even in the high wavenumber limit. As for the ranks, we observe that they slowly increase as the wavenumber grows. This relation is roughly linear (Note that the $\kappa$ axis is logarithmic), which is consistent with results previously reported in [35, 41]. The increase in ranks is to be expected as a higher wavenumber implies that less high-frequency information will be diffused by the elliptic nature of the PDE.

In the following, we examine how the preconditioner performs under $h$- and $p$-refinement. To this end, we run our experiments with various discretizations on the square domain with a fixed wavenumber of $\kappa = 19.5$. Figure 10 depicts the number of GMRES iterations and the maximal HSS rank among Schur complements under $h$-refinement. We note that the GMRES iterations increase only slightly when the problem size increases. The same holds for the HSS ranks. These results suggest that the method remains efficient as both ranks and iterations grow slowly when the problem size is increased. This is further investigated in Subsection 4.6. We conduct the same test, but perform $p$-refinement rather than $h$-refinement. Figure 11 depicts the number of GMRES iterations and the maximal HSS rank among Schur complements at the corresponding problem size $n$. For a mesh width of $h = 1/16$, we increase the polynomial degree from 3 to 12, for $h = 1/32$ from 2 to 8 and for $h = 1/64$ from 1 to 6. We notice that the growth of HSS ranks and GMRES iterations is approximately linear in the degrees of freedom.

[4] reports that the off-diagonal ranks of the inverse $A^{-1}$ scale as

$$k = \mathcal{O}(\log^{d+1} L/\epsilon_h), \tag{4.3}$$

where $L$ is the depth of the HSS cluster tree and $\epsilon_h = \mathcal{O}(h^{p+1})$ the approximation error of the Galerkin discretization. Assuming that the Schur complements exhibit similar behavior as the inverse, we expect logarithmic growth under $h$-refinement and linear growth under $p$-refinement, as confirmed in Figure 10 and Figure 11.

**4.3. Nonstandard domains.** So far, we have only considered problems formulated on the square domain $\Omega = [-1, 1]^2$. To understand how the preconditioner performs on non-standard domains, we apply it to the Helmholtz problem (4.2) formulated on a guitar shaped domain, depicted in Figure 12. The domain has a height of $H = 40$ and a width of 16. We choose an average mesh width of $h_0 = H/100$ and a polynomial degree of 4. This amounts to a total of 71760 degrees of freedom. As boundary conditions, we choose zero Dirichlet boundary conditions $g_D = 0$ on the exterior of the guitar and zero Neumann boundaries $g_N = 0$ at the sound hole. As in previous examples, we choose a constant right-hand side $f = 1$.

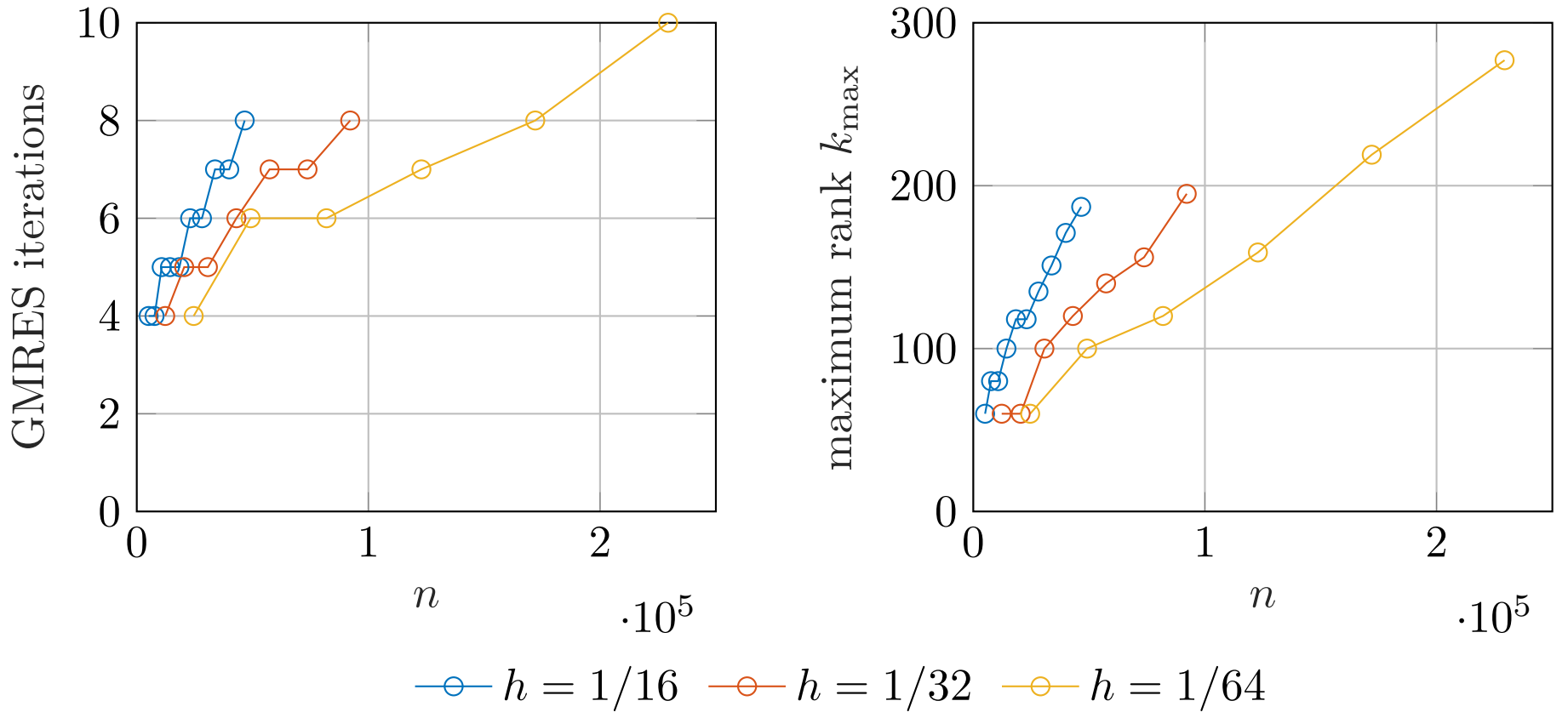

FIG. 11. *Preconditioner performance for the Helmholtz problem under p-refinement. n is the number of degrees of freedom.*

The nested dissection of the domain is depicted on the left in Figure 12, where different colors imply different nodes in the elimination tree. On the right, Figure 12 depicts the computed solutions for the wavenumbers 1, 4 and 8. We repeat previous experiments to investigate how the preconditioner performs in the high frequency limit. Figure 13 depicts the residual history of preconditioned GMRES for various wavenumbers and for various depths of the preconditioner hierarchy (the HSS switching level). This time, the effect observed in Figure 8 is more pronounced. We see that the performance deteriorates at high wavenumbers. However, this can be remedied by increasing the depth of the preconditioner hierarchy $L_{\text{HSS}}$. As we subdivide the domain, the wavenumber relative to the size of the domain becomes lower and in the limit, we recover the Poisson problem as the wavenumber approaches zero. In other words, Schur complements at lower levels can be better approximated by HSS matrices, due to the lower relative wavenumber. Consequently, we have to ensure that the hierarchy is deep enough to allow convergence in the high wavenumber limit. In our experience, it suffices to switch to HSS compression at a fixed number of elements per box, as we have done throughout the paper. In this way, the depth of the hierarchy is automatically adapted to match the quality of the discretization, which needs to match the wavenumber of the problem. A situation in which this might fail are high frequencies in combination with high polynomial orders of approximation, as our method does not account for $p$-refinement.

**4.4. Elastic wave equation.** We consider the time-independent portion of the elastic wave equations in an inhomogeneous but isotropic medium: Given a right-hand

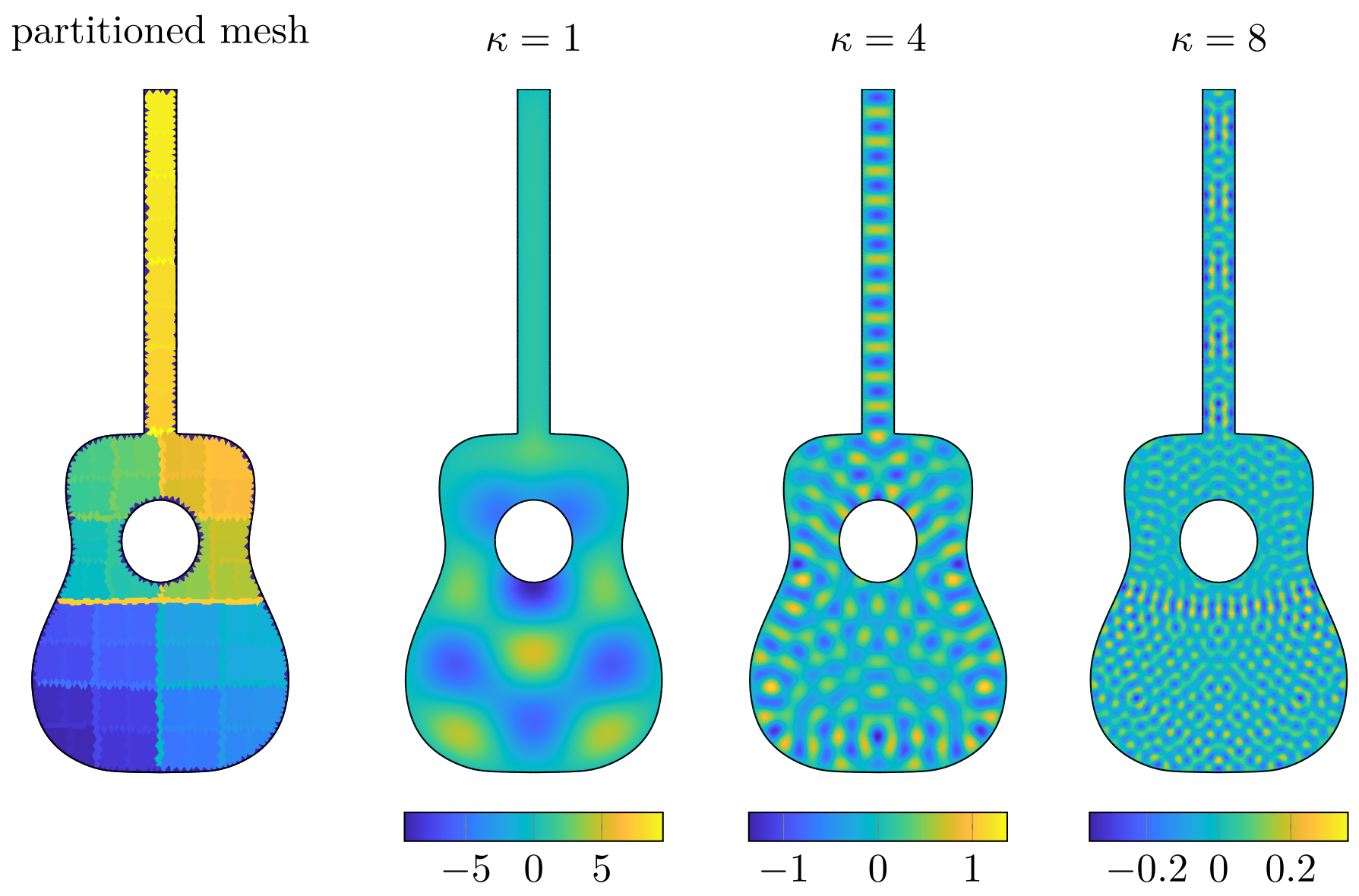

FIG. 12. *Solution of the Helmholtz problem (4.2) on a guitar shaped domain. The figure on the left depicts the nested dissection reordering. The remaining figures illustrate solutions obtained with the wave numbers* 1, 4 *and* 8.

side $f : \Omega \to \mathbb{R}^2$, find a solution $u : \Omega \to \mathbb{R}^2$ that satisfies

$$-\nabla \cdot \sigma(u) = f \text{ in } \Omega \tag{4.4a}$$

$$u = g_D \text{ on } \Gamma_D \tag{4.4b}$$

$$\sigma(u) \cdot \hat{n} = g_N \text{ on } \Gamma_N, \tag{4.4c}$$

such that the Dirichlet and Neumann boundary conditions, given by $g_D : \Gamma_D \to \mathbb{R}^2$ and $g_N : \Gamma_N \to \mathbb{R}^2$, are satisfied. Here, $\sigma(\cdot)$ denotes the two by two Cauchy stress tensor given by

$$\sigma(u) = \lambda(\nabla \cdot u)I + \mu(\nabla u + (\nabla u)^T), \tag{4.5}$$

where $\lambda, \mu : \Omega \to \mathbb{R}$ are the Lamé parameters of the material and $I$ denotes the two by two identity matrix.

Again we discretize the problem using an IPDG formulation, as presented in [8, 9]. We consider a heterogeneous problem on the square domain $\Omega = \Omega_1 \cup \Omega_2$, with $\Omega_1 = [-1, 0] \times [-1, 1]$ and $\Omega_2 = (0, 1] \times [-1, 1]$, where the material parameters are constant in each of the respective subdomains: $\mu = \mu_1$, $\lambda = \lambda_1$ in $\Omega_1$ and $\mu = \mu_2$, $\lambda = \lambda_2$ in $\Omega_2$. We choose $\mu_1 = 1$, $\mu_2 = 2$, $\lambda_1 = 1$ and $\lambda_2 = 2$. As source term we choose the constant $f = \begin{bmatrix} 0 & 1 \end{bmatrix}^T$ and we set zero Dirichlet boundary conditiones on all sides. The mesh is a regular mesh and conforms to the interface in the center of the domain.

We test the performance of the preconditioner on the described problem under $h$-refinement. Figure 14 shows that the number of GMRES iterations only increases marginally when performing $h$-refinement. On the other hand, we observe that the ranks behave roughly as in Figure 10 with the Helmholtz problem. As (4.4) is two-dimensional, the problem is double the size and we have roughly the same relative ranks as in the case of the Helmholtz problem. A similar situation presents itself under $p$-refinement, as can be seen in Figure 15. Again, the number of GMRES iterations remains roughly the same, while the ranks increase in a more or less linear manner, consistent with the behavior of the Helmholtz problem.

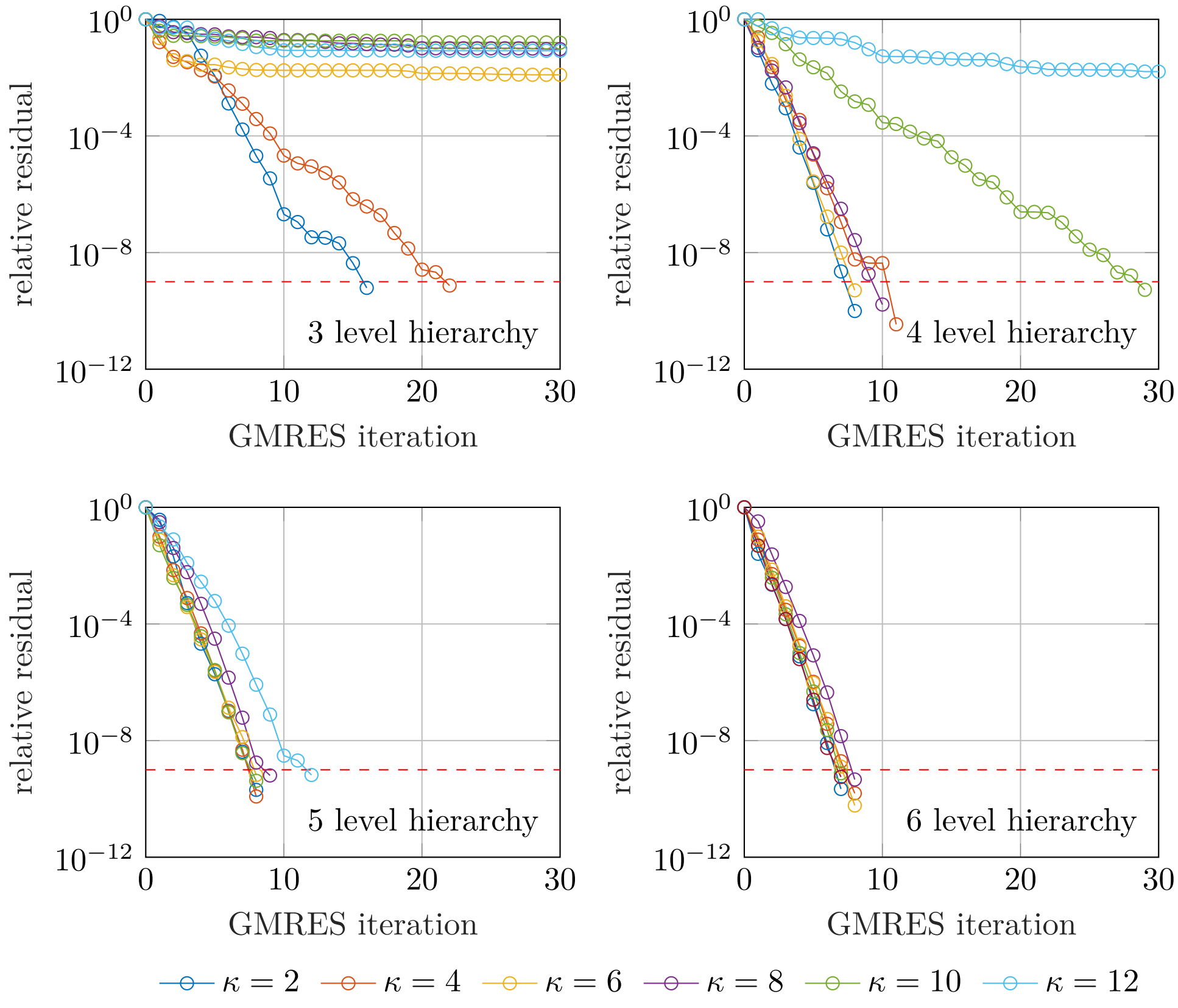

FIG. 13. *GMRES residual for various wave numbers on the guitar shaped domain with $h = 1/100$ and $p = 4$.*

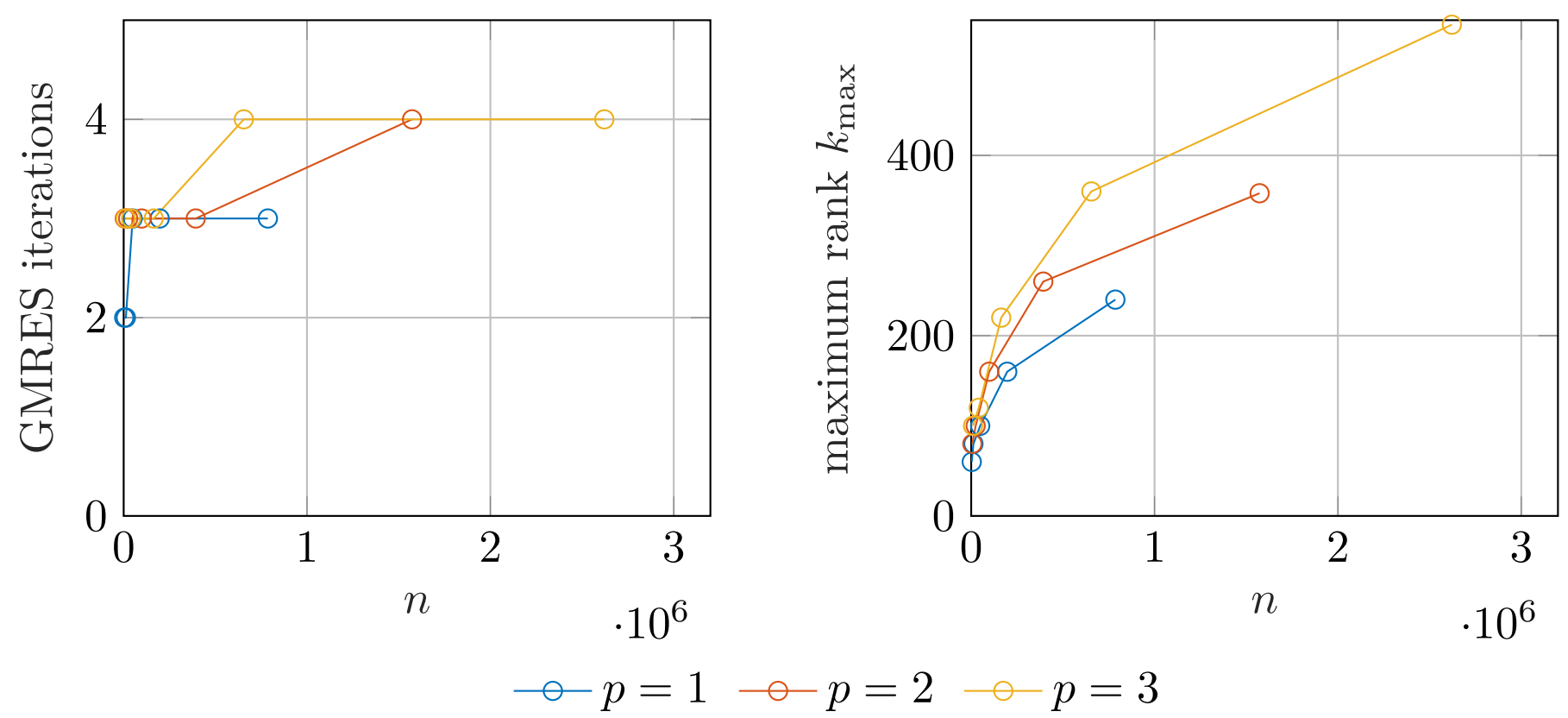

FIG. 14. *Preconditioner performance under $h$-refinement for the elastic wave equations* (4.4).

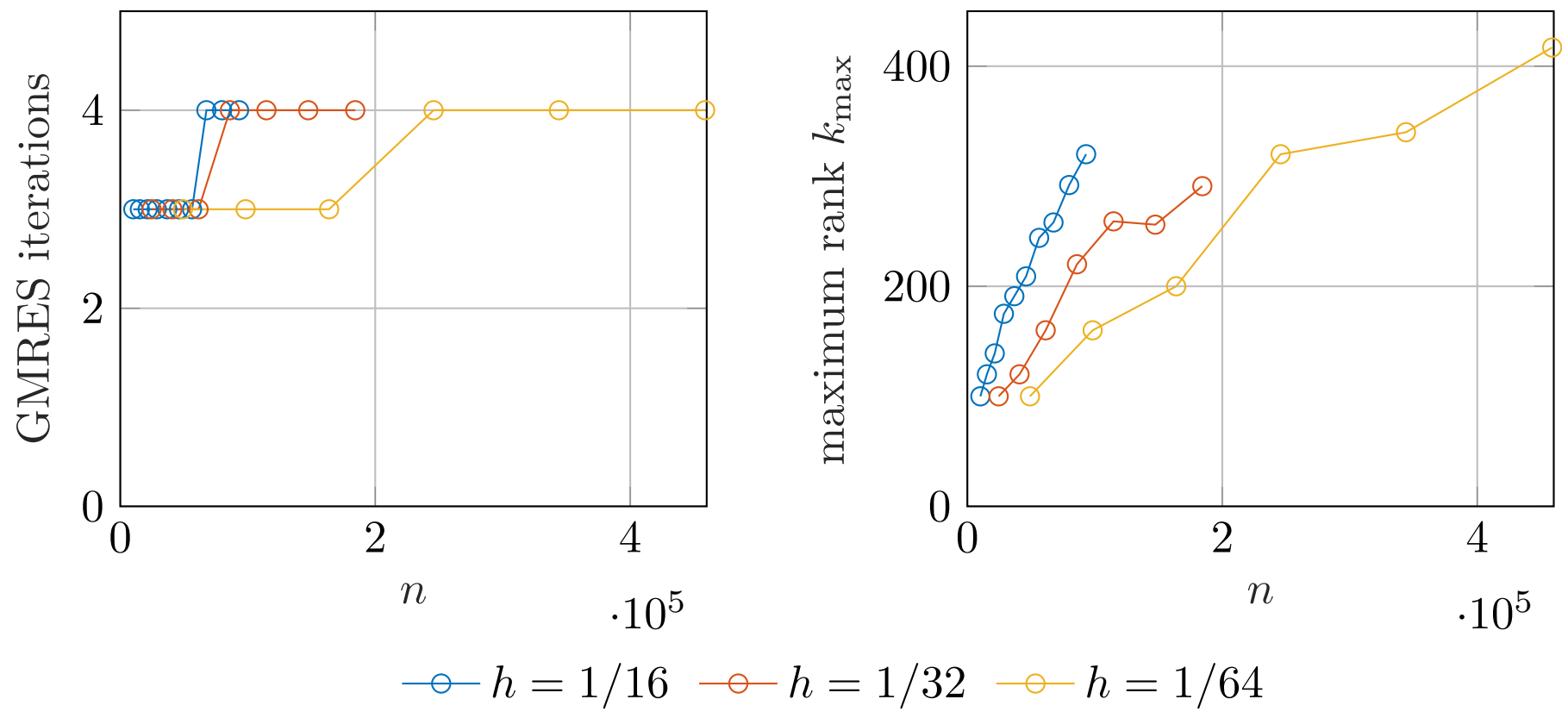

Fig. 15. *Preconditioner performance under p-refinement for the elastic wave equations* (4.4).

**4.5. Frequency-domain elastic wave equation.** As a final test, we apply the method to IPDG discretizations of the frequency-domain elastic wave equation (4.6). It can be obtained from (4.4) by applying a Fourier transform with angular frequency $\omega$. The problem reads: Find $u : \Omega \to \mathbb{R}^2$ that satisfies

$$-\nabla \cdot \sigma(u) - \rho\omega^2 u = f \text{ in } \Omega \tag{4.6a}$$

$$u = g_D \text{ on } \Gamma_D \tag{4.6b}$$

$$\sigma(u) \cdot \hat{n} = g_N \text{ on } \Gamma_N. \tag{4.6c}$$

Here, we have introduced the material density $\rho : \Omega \to \mathbb{R}$ and the angular frequency $\omega \in \mathbb{R}$. Such problems arise in the context of subsurface modeling and direct waveform inversion [28]. As such, it represents a more application-oriented test for the preconditioner.

We consider the Marmousi II velocity model [32], which is a common benchmark problem for such applications. It models soil deposits off the coast of Madagascar. The original model is 17km wide and 3.5km deep. We use the soil portion of the model which is contained in the computational domain $\Omega = [0, 17000] \times [-3500, -450]$. As boundary conditions, we choose zero Dirichlet boundary conditions $g_D = 0$ on the bottom and zero Neumann boundary conditions $g_N = 0$ on the other three sides. The source term is a dipole of the form $f = (r - r_s)\exp\big(|r - r_s|^2/(2R^2)\big)$, where $r = \begin{bmatrix} x & y \end{bmatrix}^T$ is the radius vector, $r_s = \begin{bmatrix} -1250 & 8500 \end{bmatrix}^T$ the source position and $R = 100$ the width of the dipole.

We use two meshes; one mesh is generated using a segmentation of the material distributions, to conform to the sharp interfaces in the material parameters. The other mesh is a simple, regular mesh. We use meshes with approximately 40 or 80 elements in the vertical direction, which amounts to a total of 54642 or 186868 elements in the case of the conforming mesh. The nested dissection is again generated by hierarchical subdivision, while keeping the aspect ratio of the boxes close to one.

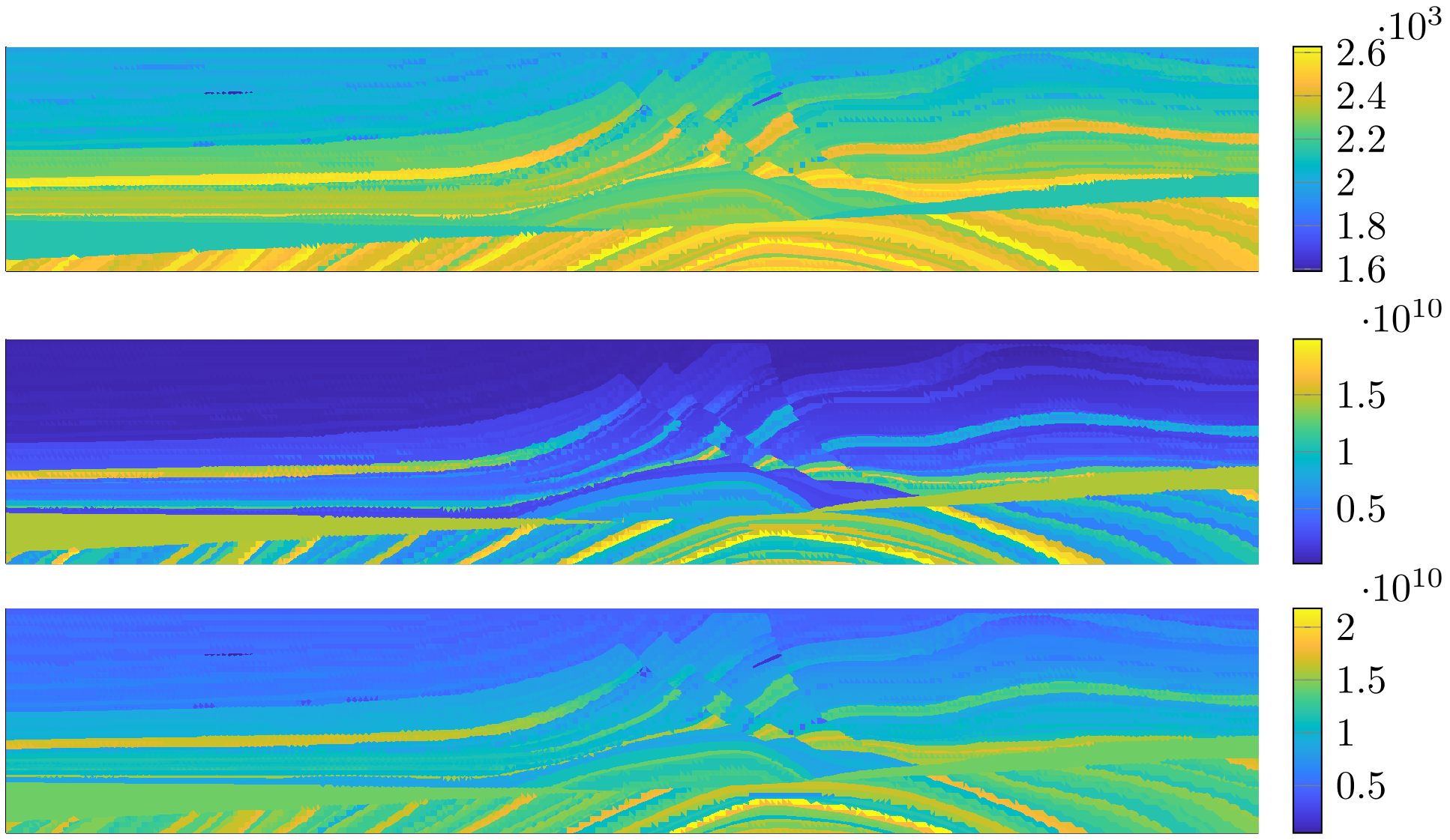

Fig. 16. *Material distribution of* $\rho$, $\mu$ *and* $\lambda$ *for* $h = 3050/40$ *and* $p = 1$.

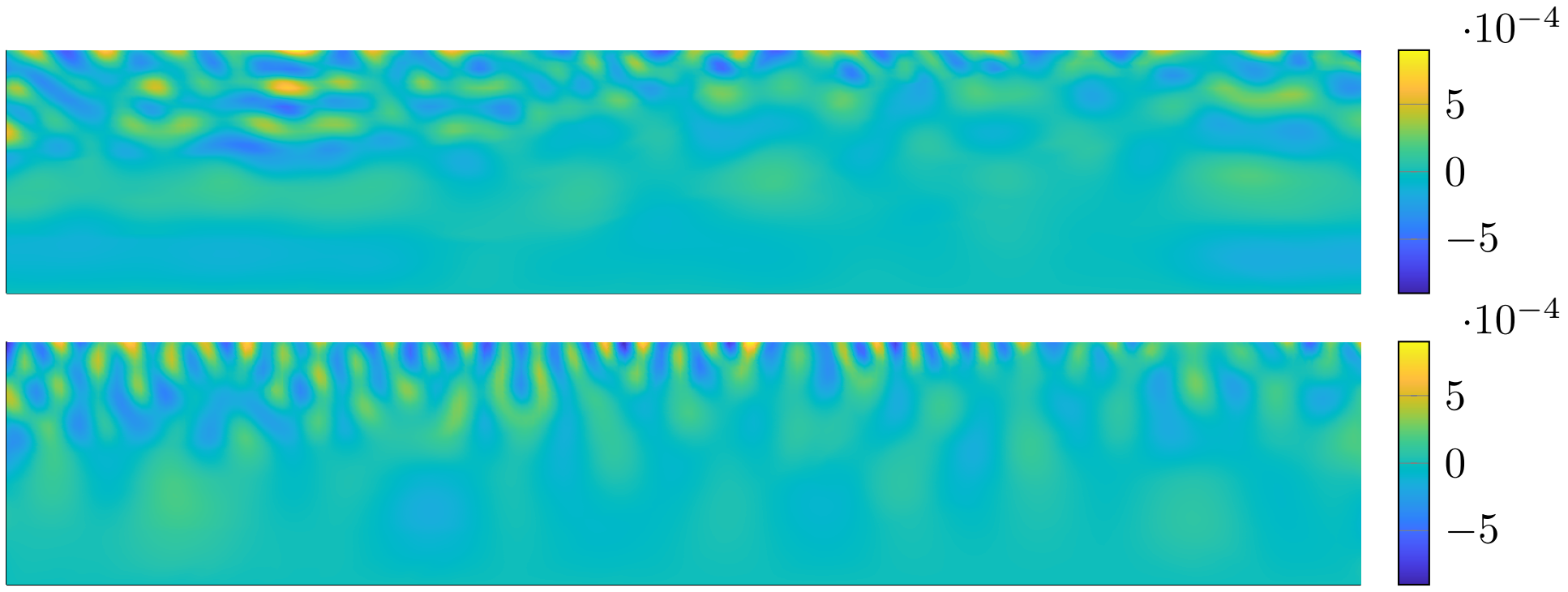

Fig. 17. $x$- *and* $y$-*components of the solution for* $h = 3050/40$ *and* $p = 1$ *and a frequency of* $\omega = 2\pi$.

For the regular mesh, the material distributions $\rho$, $\mu$ and $\lambda$ are approximated by piecewise constant functions in each element. These distributions are depicted in Figure 16. For the regular mesh, the material distributions are approximated within the discontinuous Galerkin function space on the mesh.

The problem is then solved for two frequencies, $\omega = 2\pi$ and $\omega = 8\pi$. The $x$- and $y$-components of the corresponding solution $u$ for $\omega = 2\pi$ are depicted in Figure 17. This problem is challenging as it includes many elongated high-contrast interfaces and is typically hard to precondition. Figure 18 depicts the relative residual at each GMRES iteration for both frequencies. For the purpose of comparison, we apply the incomplete $LU$ factorization (ILU) as a preconditioner, as it represents a popular, general-purpose preconditioning technique. We observe that using our method, we are able to precondition the problem and achieve a satisfying convergence rate, while the ILU preconditioner fails. This illustrates the robustness of the hierarchical pre-conditioner. The only problem-specific knowledge that is required is a hierarchical partitioning of the domain. If this partitioning is badly chosen, the method might fail or deteriorate due to Schur complements not being compressible. As we have computed the nested dissection based purely on the aspect ratio of the bounding boxes, we have demonstrated that it is quite robust, even if there are high-contrast interfaces in the domain.

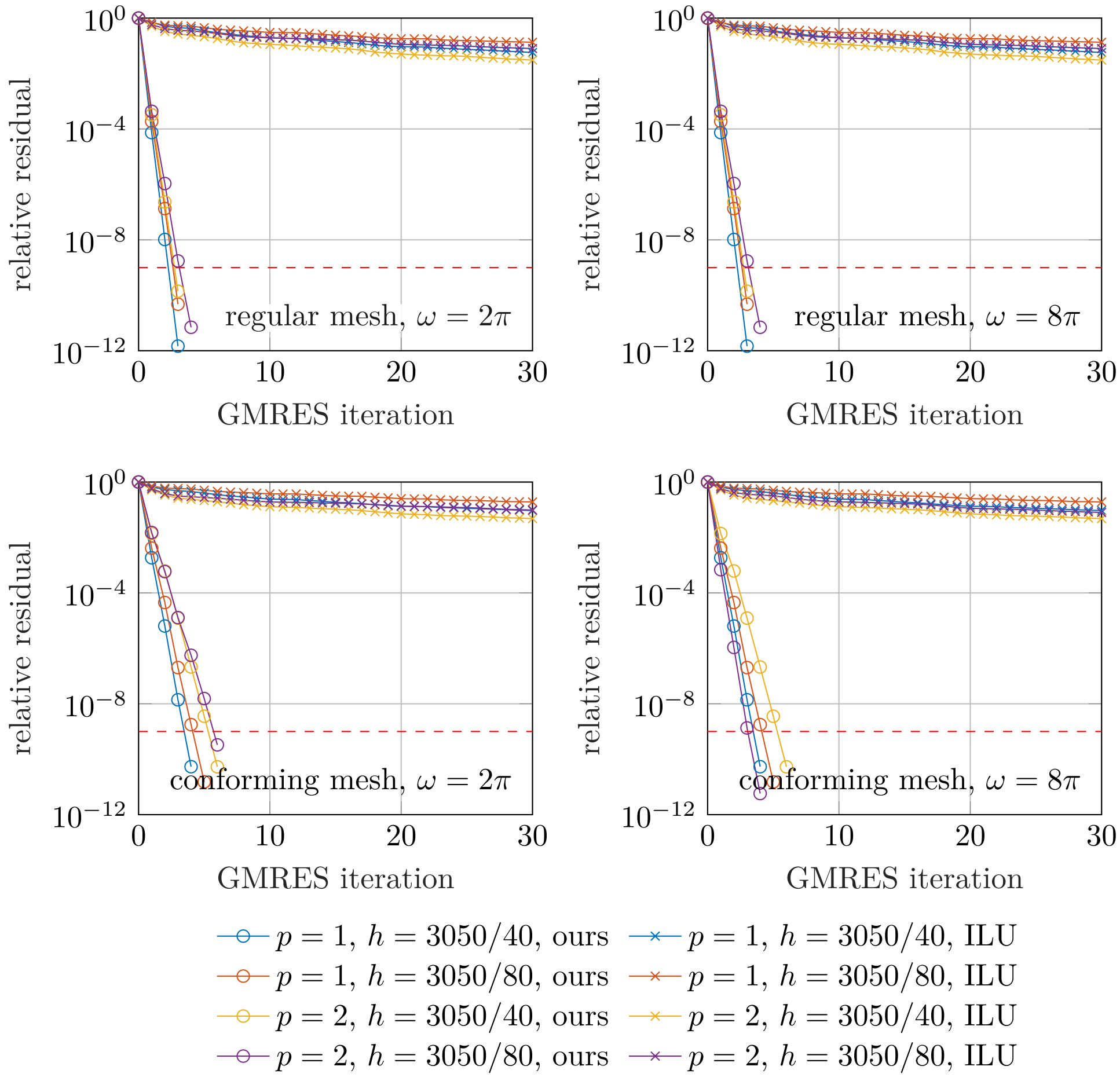

FIG. 18. *Performance of the preconditioner compared to ILU. Number of GMRES iterations for $\omega = 2\pi$ and $\omega = 8\pi$ with non-conforming and conforming discretizations.*

**4.6. Scaling.** Table 2 depicts the scaling behaviour of the preconditioner. This entails elapsed times as well as memory requirements for the construction of the approximate factorization $A \approx P = LDR$ and its application to a vector $x \to P^{-1}x$. These experiments are carried out for the two-dimensional Poisson and Helmholtz problems respectively, on meshes ranging from $h = 1/22$ to $h = 1/512$ with $p = 1$. All experiments were computed using a reasonably optimized implementation in the programming language JULIA. This data is also visualized in Figure 19.

We observe that the cost of application, construction and memory requirements roughly double from row to row, which confirms the (quasi-)linear complexity, pos-

Table 2

*Factorization and application times, as well as memory consumption for the Poisson and Helmholtz problems under h-refinement in two dimensions.*

(a) Poisson problem

| $1/h$ | $n$ | $L_{\text{HSS}}$ | application (s) | factorization (s) | memory (B) | iters | $k_{\max}$ |
|---|---|---|---|---|---|---|---|
| 22 | 2904 | 4 | 0.00311 | 0.03921 | 4.85300$e$7 | 3 | 44 |
| 32 | 6144 | 5 | 0.00684 | 0.36351 | 1.45754$e$8 | 3 | 79 |
| 45 | 12150 | 7 | 0.01442 | 0.27612 | 3.59531$e$8 | 3 | 81 |
| 64 | 24576 | 7 | 0.05801 | 0.50404 | 8.29170$e$8 | 4 | 93 |
| 90 | 48600 | 9 | 0.10054 | 1.17373 | 1.86551$e$9 | 4 | 100 |
| 128 | 98304 | 9 | 0.41032 | 2.42477 | 3.96594$e$9 | 4 | 106 |
| 181 | 196566 | 11 | 0.32823 | 5.50095 | 8.42433$e$9 | 4 | 106 |
| 256 | 393216 | 11 | 0.64399 | 18.6995 | 1.75711$e$10 | 4 | 106 |
| 362 | 786264 | 13 | 1.91412 | 40.9561 | 3.66618$e$10 | 6 | 106 |
| 512 | 1572864 | 14 | 4.00240 | 114.050 | 7.50033$e$10 | 6 | 106 |

(b) Helmholtz equation with $\kappa = 10$

| $1/h$ | $n$ | $L_{\text{HSS}}$ | application (s) | factorization (s) | memory (B) | iters | $k_{\max}$ |
|---|---|---|---|---|---|---|---|
| 22 | 2904 | 4 | 0.00318 | 0.04152 | 4.92942$e$7 | 3 | 48 |
| 32 | 6144 | 5 | 0.00744 | 0.10067 | 1.48187$e$8 | 3 | 82 |
| 45 | 12150 | 7 | 0.01418 | 0.21265 | 3.65904$e$8 | 3 | 82 |
| 64 | 24576 | 7 | 0.03072 | 0.51274 | 8.37186$e$8 | 4 | 95 |
| 90 | 48600 | 9 | 0.10588 | 1.15436 | 1.90016$e$9 | 4 | 100 |
| 128 | 98304 | 9 | 0.12941 | 2.98963 | 3.98878$e$9 | 4 | 106 |
| 181 | 196566 | 11 | 0.25293 | 5.99976 | 8.61678$e$9 | 5 | 106 |
| 256 | 393216 | 11 | 0.64557 | 19.0612 | 1.76690$e$10 | 5 | 106 |
| 362 | 786264 | 13 | 1.96164 | 42.4397 | 3.69807$e$10 | 7 | 106 |
| 512 | 1572864 | 13 | 3.69466 | 114.259 | 7.51682$e$10 | 7 | 106 |

tulated in Subsection 3.6. As the problem size increases, a bigger portion of the hierarchy is processed using compressed arithmetic using both HSS and low-rank matrices, which keeps the cost quasilinear. Moreover, we observe that the rank growth is moderate with increasing problem size which helps keep the cost low. Finally, we note a slight increase in the factorization time for the last three values of $h$. This is due to the high memory requirement, which hinders the factorization to be stored entirely in memory. As a consequence, virtual memory is used, which comes at a performance penalty.

A much harder benchmark for wave problems is to increase the wavenumber while also refining the mesh to ensure that the problems are well-resolved. As we have observed in Subsection 4.2, we expect the growth of the off-diagonal ranks to eventually increase the cost beyond the quasilinear complexity. This would also be in line with what other authors have observed with similiar methods [42, 35].

Table 3 shows timings and memory requirements for the Helmholtz equation under $h$-refinement, where we adapt $\kappa$ to keep $\kappa h$ constant. The data is also shown in Figure 20 and illustrates that the cost of factorization is roughly $\mathcal{O}(n^{4/3})$. Memory requirements appear to grow as $\mathcal{O}(n \log n)$, with $\kappa h = 1/2$ requiring more memory than $\kappa h = 1/4$, as expected.

**4.7. 3D problems.** We revisit the Poisson problem (4.1), this time in three dimensions on $\Omega = [0, 1]^3$. Problems in three dimensions are challenging for a multitude of reasons. The cost of direct methods is increased, based alone on the lower sparsity in three dimensions. For rank-based methods, it is known that the compressibility of

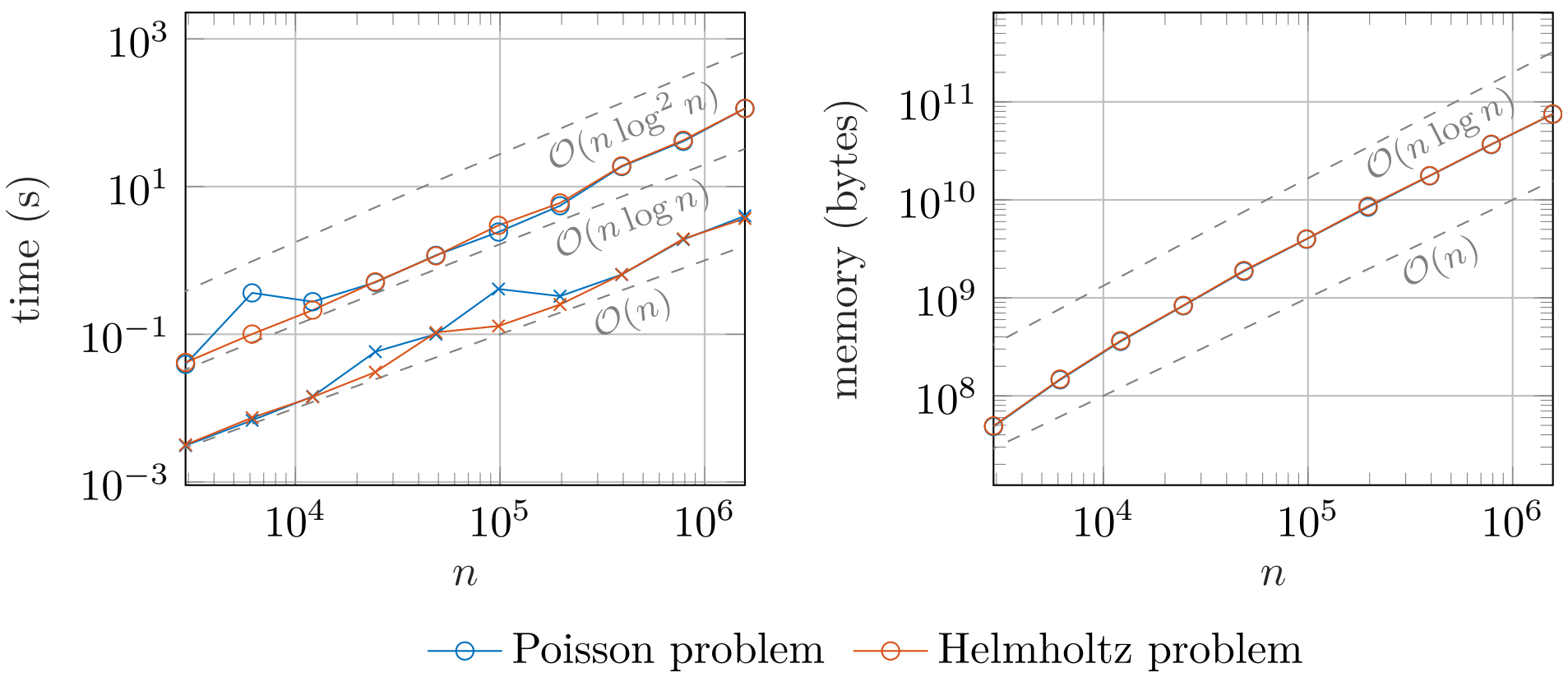

FIG. 19. *Timings and memory requirements of the preconditioner under h-refinement in two dimensions. The figure on the left shows factorization times on the top and application times on the bottom, while the figure on the right shows memory requirements.*

TABLE 3
*Factorization and application times, as well as memory consumption for the Helmholtz problem under h-refinement, while keeping $\kappa h$ constant.*

(a) Helmholtz equation with $\kappa h = 1/4$

| $1/h$ | $n$ | $L_{\text{HSS}}$ | application (s) | factorization (s) | memory (B) | iters | $k_{\max}$ |
|---|---|---|---|---|---|---|---|
| 32 | 6144 | 5 | 0.00616 | 0.06750 | 1.47068$e$8 | 3 | 80 |
| 45 | 12150 | 7 | 0.01464 | 0.27935 | 3.66537$e$8 | 4 | 81 |
| 64 | 24576 | 7 | 0.03196 | 0.50203 | 8.47778$e$8 | 3 | 94 |
| 90 | 48600 | 9 | 0.06308 | 1.20064 | 1.96278$e$9 | 4 | 105 |
| 128 | 98304 | 9 | 0.13711 | 3.12739 | 4.13341$e$9 | 4 | 106 |
| 181 | 196566 | 11 | 0.26719 | 5.46857 | 9.08292$e$9 | 5 | 106 |
| 256 | 393216 | 11 | 0.69663 | 20.3759 | 1.95125$e$10 | 5 | 146 |
| 362 | 786264 | 13 | 2.11069 | 50.1623 | 4.45347$e$10 | 7 | 186 |
| 512 | 1572864 | 13 | 4.51857 | 151.391 | 9.72291$e$10 | 6 | 266 |

(b) Helmholtz equation with $\kappa h = 1/2$

| $1/h$ | $n$ | $L_{\text{HSS}}$ | application (s) | factorization (s) | memory (B) | iters | $k_{\max}$ |
|---|---|---|---|---|---|---|---|
| 32 | 6144 | 5 | 0.00870 | 0.08765 | 1.49826$e$8 | 3 | 82 |
| 45 | 12150 | 7 | 0.01629 | 0.19147 | 3.79398$e$8 | 3 | 89 |
| 64 | 24576 | 7 | 0.03264 | 0.66125 | 8.79510$e$8 | 3 | 104 |
| 90 | 48600 | 9 | 0.07302 | 1.21613 | 2.07273$e$9 | 4 | 106 |
| 128 | 98304 | 9 | 0.14631 | 3.08951 | 4.73267$e$9 | 4 | 146 |
| 181 | 196566 | 11 | 0.39223 | 6.76212 | 1.12051$e$10 | 5 | 186 |
| 256 | 393216 | 11 | 0.76545 | 24.8398 | 2.50636$e$10 | 4 | 226 |
| 362 | 786264 | 13 | 2.29401 | 69.9011 | 6.28279$e$10 | 6 | 306 |
| 512 | 1572864 | 13 | 4.63122 | 225.629 | 1.43896$e$11 | 6 | 426 |

Schur complements decrease considerably due to the two-dimensional nature of the separators [35]. Increased ranks imply that compressed matrices have to be larger in order for the compression to be effective. This increases the overall size of problems that have to be considered to evaluate the effectiveness of our method.

We adapt the compression parameters to $\beta = 180(p+1)(p+2)(p+3)/6$ and $\epsilon_{\text{HSS}} = 10^{-2}$. Moreover, we switch to compressed arithmetic only once matrices are four times larger than the HSS blocksize $\beta$. Table 4 lists our findings. As expected, we notice that ranks are much larger, and consequently, memory quickly becomes an issue for our sequential code. The overall cost to form the approximate factorization seems to scale as $\mathcal{O}(n^2)$, whereas the cost of applying it seems to scale as $\mathcal{O}(n^{3/2})$. These findings are hardly conclusive however, due to memory limitations and the small sizes of the matrices. Nevertheless, they are consistent with previous finidings of similar rank-based approaches [35].

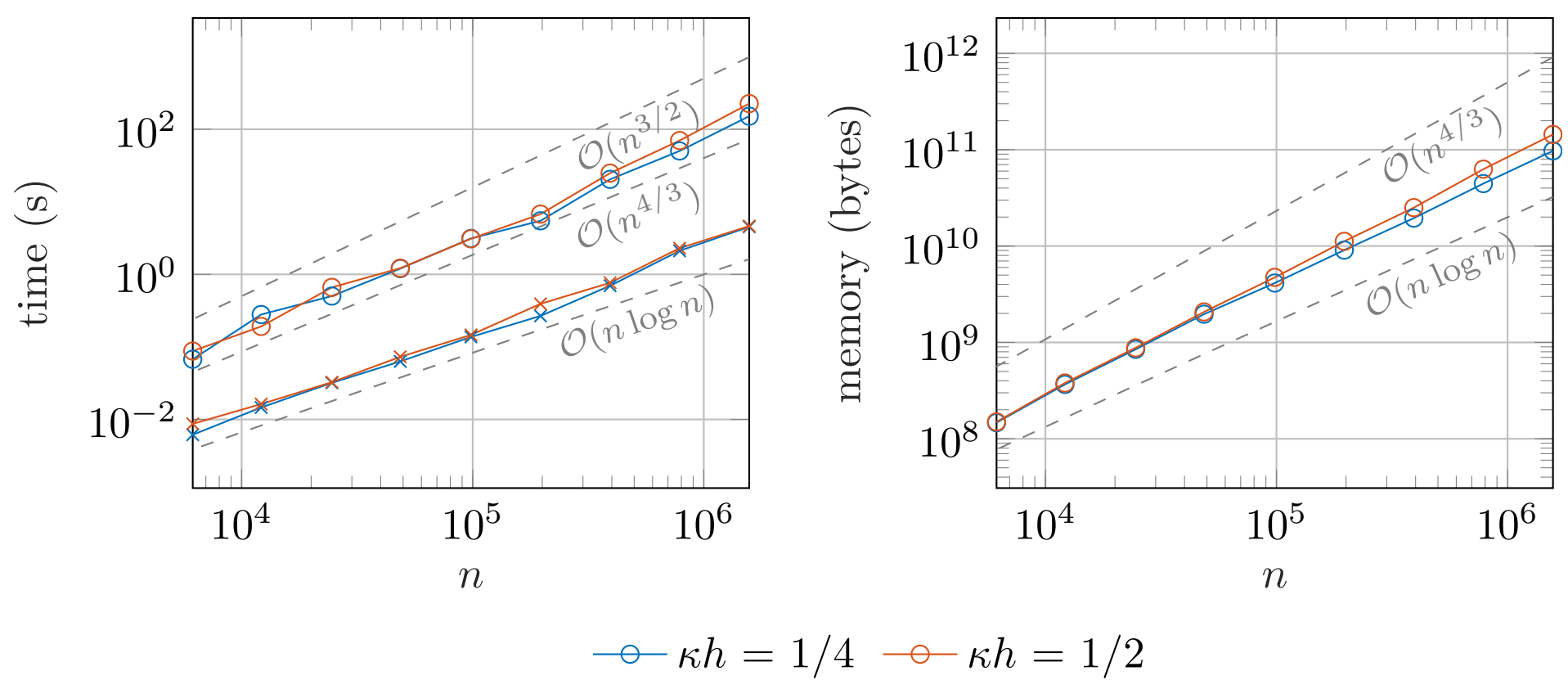

FIG. 20. *Timings and memory requirements of the preconditioner under simultaneous $\kappa$- and $h$-refinement, keeping $\kappa h$ constant.*

TABLE 4
*Factorization and application times, as well as memory consumption for the Poisson in three dimensions under $h$-refinement.*

| $1/h$ | $n$ | $L_{\text{HSS}}$ | application (s) | factorization (s) | memory (B) | iters | $k_{\max}$ |
|---|---|---|---|---|---|---|---|
| 12 | 41472 | 1 | 0.52092 | 23.6715 | 1.45695$e$10 | 10 | 780 |
| 16 | 98304 | 2 | 1.89304 | 107.020 | 4.54741$e$10 | 12 | 971 |
| 20 | 192000 | 4 | 4.41101 | 415.628 | 1.39490$e$11 | 17 | 1204 |

**5. Conclusions.** We have presented a hierarchical preconditioning technique based on nested dissection and accelerated matrix arithmetic using hierarchical matrices. This is done by exploiting properties of the finite element problem to incorporte simple HSS arithmetic into the factorization. Consequently, we can form the preconditioner matrix-free and in quasilinear $\mathcal{O}(k^2 n \log n)$ complexity, as well as apply it in $\mathcal{O}(kn)$ complexity. This approach is truly general as it only requires the original matrix $A$ and a nested dissection type elimination tree $\mathcal{E}$, which satisfies the well-separated property.

The effectiveness of this approach has been demonstrated on a variety of elliptic problems, including wave problems in the high wavenumber limit. For such problems, we observe that the depth of the hierarchy can be adapted to account for the high wavenumber. Moreover, we have shown that the preconditioner performs well in more complex situations involving heterogeneous material distributions with high contrasts and complex geometries. Consequently, we believe that the generality, efficiency and flexibility of the preconditioner make it an interesting choice for two-dimensional elliptic problems involving wave phenomena.

**Appendix A. HSS compression via random sampling.** We present a mod-ification to the randomized compression algorithm presented in [34]. The algorithm computes a probabilistic estimate of the norm [24] and performs the compression until the norm estimates fulfill a prescribed tolerance. Similar adaptive algorithms have been presented in [23]. An implementation of the algorithm can be found in [5].

Given a symmetric matrix $A$, an initial guess of the HSS rank $k$, a step size $r$ and a tolerance $\epsilon_{\text{HSS}}$, Algorithm A.1 computes an HSS approximation $A_{\text{HSS}}$ to the matrix $A$. Algorithm A.1 uses the interpolative decomposition

**Algorithm A.1** Adaptive HSS compression of a symmetric matrix $A$

Generate initial $n \times (k+r)$ random Gaussian matrix $\Omega$
Evaluate $S \leftarrow A\Omega$ via matrix-vector multiplication
**loop**
  **for all** levels $l$, starting from $L$ to 1 **do**
    **for all** nodes $i$ on level $l$ **do**
      **if** $i$ is a leaf node **then**
        $I_{\text{loc}} \leftarrow I_i^l$
        $\Omega_{\text{loc}} \leftarrow \Omega(I_i^l, :)$
        $S_{\text{loc}} \leftarrow S(I_i^l, :) - A(I_i^l, I_i^l)\Omega_{\text{loc}}$
        $D_i^{(l)} \leftarrow A(I_i^l, I_i^l)$
      **else**
        Let $c_1$ and $c_2$ be the two children of node $i$
        $I_{\text{loc}} \leftarrow [\tilde{I}_{c_1}^{l+1}, \tilde{I}_{c_2}^{l+1}]$
        $\Omega_{\text{loc}} \leftarrow \begin{bmatrix} \Omega_{c_1} \\ \Omega_{c_2} \end{bmatrix}$
        $S_{\text{loc}} \leftarrow \begin{bmatrix} S_{c_1} - A(\tilde{I}_{c_1}^l, \tilde{I}_{c_2}^l)\Omega_{c_2} \\ S_{c_2} - A(\tilde{I}_{c_2}^l, \tilde{I}_{c_1}^l)\Omega_{c_1} \end{bmatrix}$
        $B_i^{(l)} \leftarrow \begin{bmatrix} 0 & A(\tilde{I}_{c_1}^{l+1}, \tilde{I}_{c_2}^{l+1}) \\ A(\tilde{I}_{c_2}^{l+1}, \tilde{I}_{c_1}^{l+1}) & 0 \end{bmatrix}$
      **end if**
      Form interpolative Decomposition $S^* \approx S_{\text{loc}}^*(:, J_i)U_i^{(l)}$
      $\Omega_i \leftarrow \left(U_i^{(l)}\right)^* \Omega_{\text{loc}}$
      $S_i \leftarrow S_{\text{loc}}(J_i, :)$
      $\tilde{I}_i^l \leftarrow I_{\text{loc}}(J_i)$
    **end for**
  **end for**
  Generate $n \times r$ random Gaussian matrix $\tilde{\Omega}$
  Evaluate $\tilde{S} \leftarrow A\tilde{\Omega}$ via matrix-vector multiplication
  **if** $\|\tilde{S} - A_{\text{HSS}}\tilde{\Omega}\|_{\text{F}}^2 / r \leq \epsilon_{\text{HSS}}^2$ **then**
    **return** $A_{\text{HSS}}$ in HSS format, defined by $U_i^{(l)}$, $D_i^{(l)}$ and $B_i^{(l)}$
  **else**
    Increase the rank: $k \leftarrow k + r$
    Update $\Omega \leftarrow [\Omega\ \tilde{\Omega}]$ and $S \leftarrow [S\ \tilde{S}]$
  **end if**
**end loop**

$$S \approx XS(:, J), \tag{A.1}$$

to computes an approximation of $S$ that is represented via the columns $S(:, J)$. This,

in combination with the recurrence relation (2.5), allows the construction of an efficient algorithm, which only relies on the evaluation of matrix-vector products $x \to Ax$, $x \to A^*x$ and the access to individual entries $A(i,j)$.

**Acknowledgments.** We wish to thank the reviewers for their time and effort in providing us with helpful comments, which helped improve and clarify the manuscript. We thank Prof. Daniel Kressner and Dr. Stefano Massei for fruitful discussions. We thank the Swiss National Science Foundation (SNSF) for partially funding this research.